\documentclass[11pt]{article}
\usepackage{amssymb,latexsym}
\usepackage{epsfig}
\usepackage{eufrak}
\usepackage{amsmath}
\usepackage{mathrsfs}
\usepackage{color}

\newtheorem{theorem}{Theorem}
\newtheorem{lemma}{Lemma}

\newtheorem{remark}{Remark}
\newtheorem{definition}{Definition}

\newcommand{\be}{\begin{equation}}
\newcommand{\ee}{\end{equation}}
\newcommand{\bee}{\begin{eqnarray*}}
\newcommand{\eee}{\end{eqnarray*}}
\newcommand{\bel}{\begin{eqnarray}}
\newcommand{\eel}{\end{eqnarray}}
\newcommand{\bec}{\begin{cases}}
\newcommand{\eec}{\end{cases}}
\newcommand{\bem}{\begin{bmatrix}}
\newcommand{\eem}{\end{bmatrix}}

\newcommand{\bed}{\begin{description}}
\newcommand{\eed}{\end{description}}
\newcommand{\bei}{\begin{itemize}}
\newcommand{\eei}{\end{itemize}}
\newcommand{\ben}{\begin{enumerate}}
\newcommand{\een}{\end{enumerate}}

\newcommand{\beL}{\begin{lemma}}
\newcommand{\eeL}{\end{lemma}}
\newcommand{\beT}{\begin{theorem}}
\newcommand{\eeT}{\end{theorem}}

\newcommand{\bpf}{\begin{pf}}
\newcommand{\epf}{\end{pf}}

\newcommand{\pfbox}{\hfill\mbox{$\Box$}}

\newenvironment{pf}{\paragraph*{Proof{\rm.}}}{\pfbox\bigskip}

\begin{document}

\title{{\bf Fast Parallel Frequency Sweeping Algorithms for
Robust ${\cal D}$-Stability Margin
\thanks{This research was supported
in part by grants from AFOSR (F49620-94-1-0415), ARO (DAAH04-96-1-0193),
and LE{\cal Q}SF (DOD/LE{\cal Q}SF(1996-99)-04).}}}
\author{Xinjia Chen and Kemin Zhou\\
Department of Electrical and Computer Engineering\\
Louisiana State University\\
Baton Rouge, LA 70803\\
chan@ece.lsu.edu \ \ kemin@ece.lsu.edu}

\date{January 1, 2002}
\maketitle

\begin{abstract}

This paper considers the robust ${\cal D}$-stability margin problem under
polynomic structured real parametric uncertainty.
Based on the work of De Gaston and  Safonov (1988), we have developed
techniques such as, a parallel frequency sweeping strategy,
different domain splitting schemes, which significantly
reduce the computational complexity
and guarantee the convergence.

\end{abstract}

\begin{center}
\end{center}

\section{Introduction}

Robustness of control systems has been one of the central issues
 in the control community in the last two decades.
 Most of the research efforts
 have been devoted to the $\mu$ framework\cite{Bal, Ba, D, NY, ZDG} and the
 Kharitonov framework\cite{B,K,TB}.
One of the well studied robustness analysis problem is
the computation of robust stability margin under
polynomic structured real parametric uncertainty.
A number of different approaches have been proposed
 in the Kharitonov framework aimed at the nonconservative computation
 of the robust stability margin.
 Among these, we recall
 the geometric programming methods \cite{VTM},  the algorithm
 based on the Routh table \cite{SP}, and the domain splitting
 approach \cite{GS,PS, SG} based on the Zero Exclusion Condition\cite{FD}.
 In general, the algorithms in \cite{SP,VTM}
 is more efficient than the algorithm in \cite{GS}.  The main reason is that
 the algorithms in \cite{SP,VTM} are essentially based on the
 Routh-Hurwitz criterion and thus only finite conditions need to be evaluated,
 while the algorithm in \cite{GS} is based on
 the Zero Exclusion Condition and thus a frequency sweeping is essential.

\begin{figure}[htb]
\centerline{\psfig{figure=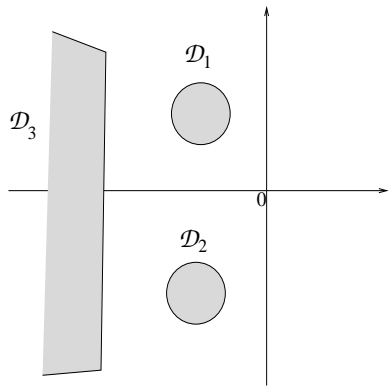,height=2in  ,width=2in  }}
\caption{Robust ${\cal D}$-Stability}
\label{fig_338}
\end{figure}

 Even though a frequency sweeping is a necessity, an algorithm based on
 the Zero Exclusion Condition has its particular advantage when dealing with
 robust ${\cal D}$-stability problems.  For example \cite{B},
for high order control systems, a typical specification might be as follows:
The closed-loop polynomials should have a pair of ``dominant roots''
in disks of given radius $\epsilon > 0$ centered
at $z_{1,2}=-u \pm j v$,
and all remaining roots having real parts less than $-\sigma$
with $\sigma > 0$ (See Figure~\ref{fig_338}, where
$z_1 \in {\cal D}_1,\;\;z_2 \in {\cal D}_2,\;\;{\cal D}=
{\cal D}_1 \bigcup {\cal D}_2 \bigcup {\cal D}_3$).
Then, a robust ${\cal D}$-stability margin problem can
be defined as follows:
 What is the maximum perturbation of plant parameters
 such that the roots of
the closed loop polynomial
remain robustly in
${\cal D}=\{z \in {\bf C}:\; |z-z_1| < \epsilon\} \bigcup
\{z \in {\bf C}:\; |z-z_2| < \epsilon\} \bigcup
\{z \in {\bf C}:\; \Re(z) < \sigma\}$?
  Since the root region ${\cal D}$ can be
defined as a union of disjoint open subsets with complicated
boundary in the complex plane, the robust ${\cal D}$-stability problems in general cannot be solved by
existing results in the $\mu$ framework or the algorithms in \cite{SP,VTM}
which are based on the Routh-Hurwitz criterion.
 For special cases that  ${\cal D}$
is simply connected and is defined via the Nyquist curve of certain
rational polynomials $f(s)=\frac{g(s)}{h(s)}$,
the robust ${\cal D}$-stability problem of $p(s)$ may be reduced to the
robust stability problem of polynomial
$\hat{p} (s) = p(f(s)).(h(s))^{n_h}$ where $n_h$ is
the degree of polynomial $h(s)$ and then the algorithms
in \cite{SP,VTM} may be applied.  However, the complexity is increased
 substantially because the coefficients of $\hat{p} (s)$ may be very complex and
 the degree of $\hat{p} (s)$
 is $n_g$ times of the degree of $p(s)$ where $n_g$ is
the degree of polynomial $g(s)$ \cite{S}.

 The advantage of an algorithm based on the Zero Exclusion Condition is
 that it can be applied to the robust ${\cal D}$-stability problem with
 arbitrary complicated root region ${\cal D}$.
 What we only need to do is to verify whether the Zero Exclusion Condition
 is satisfied for all boundary point of ${\cal D}$.
 In situations where the root region ${\cal D}$ is complicated,
 applying such algorithms becomes essential.
 However, the computational complexity can be very high.
 In particular, the growth of computations is exponential with the number of parameters.
 It has been shown that these type of problems are in general NP-hard (see \cite{Rohn} and the reference therein).
Moreover, since the frequency search can be discontinuous (as shown in {\cite{BKT}),
we usually need to evaluate the Zero Exclusion Condition
 for many boundary points of ${\cal D}$ to come up with a
 reasonably accurate solution.  Therefore, there is strong motivation
 to develop efficient algorithms based on the Zero Exclusion Condition to
 tackle the robust ${\cal D}$-stability problems.

The algorithm proposed by De Gaston and  Safonov \cite{GS} is based on the
Zero Exclusion Condition and thus can be applied to the general
robust ${\cal D}$-stability problem.  However, there exist two problems.

First of all, it is noted that the convergence of the
algorithm in \cite{GS} was concluded upon an impractical assumption.
 That is,  a domain can be divided fine enough to converge to a point
(see \cite{GS} line $40-53$ of page $156$ in the proof of the
Convergence Theorem).
However, to satisfy this assumption,  the computation complexity may be
unacceptably high.  In this paper, we have shown that it is sufficient to
guarantee the convergence in computing the stability margin by
guaranteeing that the distance between critical vertices converge to $0$.
Therefore, it is not necessary to divide a subdomain so many times to
collapse it to a point.  In contrast, what we need is to make the critical
vertices crunch together.  Thus, the computation can be reduced greatly.
We provide two splitting schemes which guarantee this.

Another problem with the algorithm in \cite{GS} is its inefficiency.
One main hurdle is its tedious frequency sweeping.  Consider a family of
uncertain polynomials $p(s,q),\;q \in {\cal Q}$ where ${\cal Q}$ is the set of
uncertain parameters.  Let $k_{m}(\omega,{\cal Q})=\sup
\{k: 0 \notin p(j \omega,k{\cal Q})\}$ where $p(j \omega,k{\cal Q})$ is the
{\it value set} associated with frequency $\omega$ and perturbation bound $k$.
 The algorithms in
\cite{GS} compute $k_{m}(\omega,{\cal Q})$ exactly for each frequency
$\omega$ and compare to find the minimum as the stability margin.
 To the best of our knowledge, all frequency sweeping
 techniques in the literature
 follow this format.  In this paper, we investigate a smart frequency
sweeping strategy.  We compute $k_{m}$ for $n_r > 1$ frequencies in parallel.
Domain splitting is also performed in  parallel at each iteration level.
 Information is exchanged among all subdomains to determine
 which subdomain for which frequency
 should be eliminated from further consideration without obtaining the exact
 value of $k_{m}$.  The stability margin is achieved as
 the minimum record of the upper bounds of all the subdomains ever generated.
 The convergence rate is much faster than that of \cite{GS}.

The paper is organized as follows. Section $2$ introduces the
robust ${\cal D}$-stability problem and
the work of De Gaston and  Safonov \cite{GS}.
Section $3$ discusses the Convergence Theorem of \cite{GS} and
different domain splitting schemes.
Section $4$ presents our Parallel Frequency Sweeping Algorithm.
An illustrative example is given in Section $5$ and
Section $6$ is the conclusion.


\section{Preliminary}

It is well known that the stability problem of an MIMO system can be reduced to
the study of the root location of a related polynomial \cite{B,C}.
We consider a family of polynomials $p(s,q)$ of degree $n$ whose coefficients
$a_{i}(q)$ are continuous functions of $\ell$-dimensional
vector of real uncertain parameters $q$, each bounded
in the interval $[q_{i}^{-}, q_{i}^{+}]$.  More formally, we define
\[
p(s,q):=a_{0}(q)+a_{1}(q) s+a_{2}(q) s^{2} + \cdots + a_{n}(q) s^{n}
\]
where $q:=(q_{1}, \cdots, q_{\ell})$ and
the hypercube
 ${\cal Q}:= \{q:\;q_{i}^{-} \leq q_{i} \leq  q_{i}^{+}, i=1,\cdots,\ell\}$
 with the nominal parameter $q^0 \in {\cal Q}$.

\subsection{Robust ${\cal D}$-stability Margin}

\begin{definition} \label{d5}
Let ${\cal D}$ be an open region in the complex plane and take $p(s)$
to be a fixed polynomial. Then $p(s)$ is said to be ${\cal D}$-stable
if all its roots lie in the region ${\cal D}$.
\end{definition}

\begin{definition} \label{d6}

A family of polynomials ${\cal P}=\{p(.,q): q \in {\cal Q}\}$
is said to be robustly
${\cal D}$-stable if all roots of $p(s,q)$ lies in ${\cal D}$.  For special case when
${\cal D}$ is
the open left half plane, ${\cal P}$ is simply said to be robustly stable.

\end{definition}

Let $Q \subseteq {\cal Q}$. Define {\it value set} $p(z, Q) \subset
{\bf C}$ as follows:
\[
p(z, Q):= \{p(z, q): q \in Q \}.
\]
Define
\[
kQ:= \{q^0 + k(q-q^0):q \in Q \}.
\]

We first state the Zero Exclusion Condition for uncertain polynomials.

\begin{theorem} \label{t1} {\bf (\cite{FD})}
The polynomial $p(s,q)$ is robustly ${\cal D}$-stable for all
$q \in {\cal Q}$ if and only if $p(s,q)$ is stable for some
$q \in {\cal Q}$ and $0 \notin p(z, {\cal Q})$ for all
$z \in \partial {\cal D}$.

\end{theorem}

Let ${\cal D}_{1},\;{\cal D}_{2},\cdots,{\cal D}_{N}$ be disjoint open
subsets of the complex plane and suppose
${\cal P}=\{p(.,q):q \in {\cal Q}\}$ is a family of polynomials with invariant
degree.  For each $q \in {\cal Q}$ and $i \in \{1,\cdots,N\}$,
let $n_{i}(q)$ denote the number of roots of $p(s,q)$ in ${\cal D}_{i}$.
Finally, assume that $p(s,q^{0})$
has no roots in the boundary of
${\cal D}=
{\cal D}_{1} \bigcup {\cal D}_{2} \bigcup \cdots {\cal D}_{N}$.
Then each of the root indices $n_i(q)$ remains invariant over ${\cal Q}$
if and only if
the Zero Exclusion Condition $0 \notin p(z,{\cal Q})$ is satisfied
for all points in $\partial {\cal D}$.

\begin{definition} \label{d7}
Suppose ${\cal D}$ is an open subset of
the complex plane with boundary ${\partial {\cal D}}$.
Then, given an interval $I \subseteq {\bf R}$,
a mapping $\Phi_{\cal D}: I \rightarrow  {\partial {\cal D}}$ is
said to be a boundary sweeping function for ${\cal D}$ if $\Phi_{\cal D}$
is continuous and onto; i.e., $\Phi_{\cal D}$ is continuous and for each
point $z \in \partial {\cal D}$, there is
some $\delta \in I$ such that $\Phi_{\cal D}(\delta) = z$.
The scalar $\delta$ is called a generalized frequency variable for ${\cal D}$.

\end{definition}

Let $k_{m}(\delta,Q):=\sup\{k:0 \notin p(\Phi_{\cal D}(\delta), k Q)\}$.
 The robust ${\cal D}$-stability
margin $k_{max}$ is given by $k_{max}=\inf_{\delta \in I}
 k_{m}(\delta,{\cal Q})$.
In general, when ${\cal D} = \bigcup _{l}^{N} {\cal D}_{l}$ where
  ${\cal D}_{l}, \;l=1,\cdots,N$  are disjointed open
  subsets in the complex plane, we can define
  $N$ boundary sweeping functions
  $\Phi_{\cal D}^{l}:I_{l} \rightarrow
  \partial {\cal D}_{l}, \;\;l=1,\cdots,N$ respectively.
  Then the robust ${\cal D}$-stability margin is given by
  \[
k_{max}=\min_{l=1,\cdots, N} \inf_{\delta \in I_{l}}  k_{m}(\delta,{\cal Q}).
\]

\subsection{Domain Splitting Algorithms}

It is noted that the analysis of robustness under
polynomic structured real parametric uncertainty
can be converted into a simpler analysis problem
dealing with multilinear structured uncertainty \cite{SG,PS}.

\begin{definition} \label{d8}
An uncertain polynomial $p(s,q)=\sum_{i=0}^{n} a_{i}(q)s^i$ is
said to have a
multi-linear uncertainty structure if each of the coefficient
functions $a_{i}(q)$ is multi-linear.  That is, if all but one
component of the
vector is fixed, then $a_{i}(q)$ is affine linear
in the remaining component of $q$.  More generally, $p(s,q)$ is said to
have a polynomic uncertainty structure if each of the coefficient
functions $a_{i}(q)$ is a multi-variable polynomial in the components of $q$.

\end{definition}

In general, there exists no analytic solution for computing exactly
$k_{max}$.  However,  the following Mapping Theorem can be applied to obtain
a lower bound for $k_m(\delta,Q)$ for a family of polynomials of
multi-linear uncertainty structure.

\begin{theorem} \label{t2} {\bf (\cite{ZD})}
Suppose an uncertain polynomial $p(s,q)$ has a
multi-linear uncertainty structure.  Then
\[
{\rm conv} \; p(z, Q) =
{\rm conv} \; \{ p( z, q^{1}),
p( z, q^{2}), \cdots,  p( z, q^{2^\ell}) \},\;\;\;\forall z \in \partial {\cal D}
\]
where conv denotes the convex hull and $q^1,\cdots,q^{2^\ell}$
denotes the $2^\ell$
vertices of the hypercube $Q$.
\end{theorem}

Let $k_{l}(\delta,Q):=\min\{k:\;0 \in {\rm conv}\;p(\Phi_{\cal D}(\delta), kQ)\}$.
Then by the Mapping Theorem, $k_{l}$ is a lower bound,
i.e., $k_{l} \leq k_{m}$.

\begin{definition} \label{d1} {\bf (\cite{GS})}
Critical vertices are those adjacent extreme points $M_\alpha, M_\beta$ of
${\rm conv} \; p(z, k_{l}Q)$ such that $0 \in conv \; \{M_\alpha, M_\beta\}$.
\end{definition}

\begin{definition} \label{d2} {\bf (\cite{GS})}
$m(\alpha,\beta)$ is the number of differing coordinates of
two vertices $q^\alpha, q^\beta$ that are mapped by
$p(z, .)$ to $M_\alpha, M_\beta$, respectively.

\end{definition}

Obviously that it follows from the Mapping Theorem that
 $k_{l}=k_{m}$ for $m(\alpha,\beta)=0,\;1$.  For $m(\alpha,\beta) \geq 2$,
 a {\it vertex path} is defined as follows.

\begin{definition} \label{d3} {\bf (\cite{GS})}

A vertex path is any path between critical vertices $M_\alpha, M_\beta$
consisting of $m(\alpha,\beta)$ straight-line segments defined by $p(z, .)$
as $q$ progresses from $q^\alpha$ to $q^\beta$ along the edges
of the hypercube $Q$.

\end{definition}

Define
\[
k_{u}(\delta,Q):=\inf \{k: {\rm At \;least \;one \;of \;the \;vertex \;
paths \;of } \;{\rm conv} \; p(\Phi_{\cal D}(\delta), kQ) \;
{\rm intercepts \; the \; origin}\}.
\]
It is shown in \cite{GS} that $k_{u}$ is an upper bound,
i.e., $k_{m} \leq k_{u}$.

In general, it is impractical to compute $k_{max}$ over all frequencies.
The techniques
developed in \cite{GS, PS, SG} work essentially as follows.

Choose a range of frequency $[\delta_l,\delta_u] \subset I$
and grid it as
\be
\delta_j:= \delta_l + \frac{(\delta_u-\delta_l)(j-1)}{n_r n_c},\;\;j=1,\cdots,n_r n_c
\label{grid}
\ee
where $n_r \geq 2,\;\;n_c \geq 1$ are integers.
Apply Algorithm $1$ to compute an upper bound $k_u^j$
and a lower bound $k_l^j$ for $k_{m}(\delta_j,{\cal Q})$ such that
$\frac{k_u^j-k_l^j}{k_l^j} < \epsilon, \;\;j=1,\cdots, n_r n_c$.
 Then an estimate of $k_{max}$ can be defined as
\[
\tilde{k}_{max}:= \min_{j=1,\cdots, n_r n_c} k_{m}(\delta_j,{\cal Q})
\]
which satisfies
\[
\min_{j=1,\cdots, n_r n_c} k_l^j
\leq \tilde{k}_{max} \leq \max_{j=1,\cdots, n_r n_c} k_u^j
\]
with
\[
\frac{ \max_{j=1,\cdots, n_r n_c} k_u^j - \min_{j=1,\cdots, n_r n_c} k_l^j }
{\min_{j=1,\cdots, n_r n_c} k_l^j } < \epsilon.
\]

{\bf Algorithm $1$ \cite{GS}--- Computing $k_{m}(\delta,{\cal Q})$ }

\begin{itemize}

\item Step $1$: Determine lower bound on $k_{m}$.  Designate the
initial uncertain parameter domain, the $n$-dimensional
hypercube ${\cal Q}$, as $Q_{11}$.

\item Step $2$: Determine upper bound on $k_{m}$.

\item Step $3$: Iterate to converge lower and upper bounds to $k_{m}$.
Establish an iterative procedure with counter $r=1,2,3,\cdots$.
For each iteration perform the following operations on
subdomains $Q_{rp}$ where $p$ represents
the number of subdomains left in consideration after the $r$th iteration.

\item Step $3$--$1$: Increment $r$, i.e., $r \leftarrow r+1$.

\item Step $3$--$2$: Make orthogonal cuts midway on the longer edges of
each subdomain $Q_{rw},\;w=1,\cdots,p$ in order that all edge
length ratios remain within a factor of $2$ of each other.
Designate these two subdomains as $Q_{rw}$ and $Q_{r(w+p)}$.

\item Step $3$--$3$: Obtain $k_{l_{rw}}:=k_l(\delta, Q_{rw})$ and
$k_{l_{r(w+p)}}:=k_l(\delta,Q_{r(w+p)})$ via Step $1$.
(Note: See \cite{GS} for handling exceptions).

\item Step $3$--$4$: Obtain $k_{u_{rw}}:=k_u(\delta, Q_{rw})$ and
$k_{u_{r(w+p)}}:=k_u(\delta,Q_{r(w+p)})$ via Step $1$.
(Note: See \cite{GS} for handling exceptions).

\item Step $3$--$5$: Repeat Steps $3$--$2$ to
$3$--$4$ for each $w=1,\cdots,p$.

\item Step $3$--$6$: Define
$
k_{l_{r}}:=
\min \{ k_{l_{r1}}, k_{l_{r2}}, \cdots, k_{l_{r(2p)}} \}
$
and
$
k_{u_{r}}:=
\min \{ k_{u_{r1}}, k_{u_{r2}}, \cdots, k_{u_{r(2p)}},\;k_{u_{r-1}} \}$ and
$\epsilon_{r}=\frac{k_{u_r}-k_{l_r}}{k_{l_r}}$.
(Note: It is shown in \cite{GS} that
$
k_{l_{r-1}} < k_{l_{r}} < k_{m} \leq k_{u_{r}} \leq k_{u_{r-1}}$.)

\item Step $3$--$7$: Eliminate from further consideration all subdomains
$Q_{rw},\;w=1,\cdots,2p$, whose associated $k_{l_{rw}} > k_{u_{r}}$.
Designate the number so eliminated as $u$ and define a new $p=2p-u$.

\item Step $3$--$8$: Repeat Steps $3$--$1$ to $3$--$8$ until
$k_{l_{r}} \rightarrow k_{m}$.
The stop criteria is that $\epsilon_{r}$ is less than a chosen
tolerance $\epsilon > 0$.

\end{itemize}

\begin{remark}
In the above conventional frequency sweeping algorithm,
the most important mechanism which impacts the efficiency is the
elimination of subdomains whose lower bounds are greater than the minimum
record of the upper bounds of all subdomains of the frequency being evaluated.
This mechanism is implemented in Step $3$--$7$.  It has been demonstrated in \cite{Pena}
that, although the theoretical increase in subdomains should be exponential,
in practice the growth can be linear due to such mechanism.
We would like to note that the elimination processes for different
frequencies are independent and hence such independent feature leaves room for
a substantial reduction of the growth of subdomains.
\end{remark}

\section{Domain Splitting and Convergence}

One of the most important requirement of an algorithm is on its convergence.
For example, for the above algorithm, it is expected that
 given any tolerance $\epsilon > 0$, the above algorithm
 stops at finite iteration, i.e., $r < \infty$ for each frequency.
 In this section,
 we investigate how a domain splitting can affect the convergence.

\subsection{An Impractical Assumption}

The convergence of the above algorithm was addressed in \cite{GS,G}
and a {\it Convergence Theorem}
was proposed.  However, in the proof of the {\it Convergence Theorem}\cite{GS},
the convergence was concluded upon the assumption that
each subdomain converges to
a single point by subdivisions (see \cite{GS},
lines $40-53$ of page $162$).  In another paper \cite{SG},
the convergence was also concluded by assuming that a subdomain
is divided fine enough (see the last paragraph in page $767$ of \cite{SG}).
In fact,  such an assumption is in general impractical to be satisfied.
This is because the computation is usually very high to divide
a domain fine enough to collapse it to a point.

In the above algorithm, the criterion adopted in splitting a domain
is that ``make orthogonal cut midway on the larger edges''.
It is also addressed in \cite{G,SG} that a splitting of a domain should be
made in a way guaranteeing two critical vertices remained in
different subdomains.  In general, there is more than one
way to satisfy these two criteria.  We would like to note that,
in general, a splitting scheme which just consists of these two criteria
is not sufficient to obtain a sequence of
lower bounds (or upper bounds) converging to $k_m$ or a sequence of
subdomains converging
to a single point in ${\bf R}^\ell$.

For example, consider a hypercube
${\cal Q}=\{q \in {\bf R}^{5} :\; q_i \in [0, 1],\; i=1,\cdots,5\}$.
Let $Q_1={\cal Q}$.  Based on the above two criteria,
$Q_{r}$ can be splitted as $Q_{r+1}$ and $Q_{r+1}^{'}$
with $k_{m}(\delta,Q_{r+1})=k_{m}(\delta,{\cal Q})$ at the $r$-th splitting,
$r=1,2,\cdots$.  We cannot exclude the possibility that
there exists $r_{c} < \infty$ such that the following are true.

\begin{itemize}

\item
Critical vertices differ in coordinates $q_{1},\; q_{2},\;q_{3},\;q_{4}$ for $ r > r_{c}$.

\item Coordinates $q_{1},\; q_{2}$ are cut in round robin order for $r > r_{c}$.

\end{itemize}

Finally, we will end up with a degenerate hypercube,
with $q_{1},\; q_{2},\;q_{5}$ being constants and $q_{3},\;q_{4}$ varying
within intervals,  i.e., a planar ``box''.  Because $q_{3},\;q_{4}$
can vary in intervals,  it is possible that there is a gap between the
upper bound $k_{u}$ and the lower bound $k_{l}$, i.e., $\exists \nu > 0$ such
that $k_{u}- k_{l} > \nu$.

Therefore, it is important to raise the following question:

{\it What kind of splitting guarantees the convergence?}

\subsection{Guaranteed Critical Vertices Distance Convergence}

Consider a hypercube
${\cal Q}=\{q \in {\bf R}^{\ell} :\; q_{i} \in
[q_i^{-}, q_i^{+}],\; i=1,\cdots,\ell\}$.
Define a sequence (finite or infinite) of domains $\{Q_r\}$
iteratively as follows.

\begin{itemize}

\item Step $1$: Let $Q_{1}={\cal Q}$. Let $r=1$.

\item Step $2$: If the critical vertices of domain $Q_r$, denoted by
$q^{\alpha_{r}}$ and $q^{\beta_{r}}$, differ in no more than one coordinate
then the iteration process is terminated, otherwise choose
$i_\star \in {\cal U}_{r}:=\{i:\;q^{\alpha_{r}}_i \neq q^{\beta_{r}}_i \}$
and designate either
$
\{q \in Q_r: \frac{q^{\alpha_{r}}_{i_\star,r} + q^{\beta_{r}}_{i_\star,r}}{2}
\leq q_{i_\star} \leq
\max \{q^{\alpha_{r}}_{i_\star,r} , q^{\beta_{r}}_{i_\star,r}\} \}
$
or
$
\{q \in Q_r:
\min \{q^{\alpha_{r}}_{i_\star,r} , q^{\beta_{r}}_{i_\star,r}\}
 \leq q_{i_\star} \leq
\frac{q^{\alpha_{r}}_{i_\star,r} + q^{\beta_{r}}_{i_\star,r}}{2}\}
$
as $Q_{r+1}$.
\item Step $3$: Set $r=r+1$. Go to Step $2$.

\end{itemize}

In general, there are more than one way to choose $i_\star \in {\cal U}_{r}$.
Let $Q_r= \{q \in {\bf R}^{\ell} :\;
q_{i} \in [q_{i,r}^{-}, q_{i,r}^{+}],\; i=1,\cdots,\ell\}$, we can
define a splitting scheme as follows.

\begin{definition} \label{d10}

A maximal-cut is a partition of $Q_r$ as above by choosing
 $i_\star \in {\cal U}_{r}$ such that
\[
{ q^+_{i_\star,r}-q^-_{i_\star,r} }
=\max _{i \in {\cal U}_{r} }
\;{ q^+_{i,r} - q^-_{i,r} }.
\]

\end{definition}
Another splitting scheme adopted in \cite{G} was
 that the cut should be made over the coordinate that has been
 subdivided the least number of times.
More formally, we define a {\it fair-cut} scheme as follows.
\begin{definition} \label{d11}
A fair-cut is a partition of $Q_r$ as above by choosing
$i_\star \in {\cal U}_{r}$ such that
\[
\frac { q^+_{i_\star}-q^-_{i_\star} }
{ q^+_{i_\star,r}-q^-_{i_\star,r} }
=\min _{i \in {\cal U}_{r} }
\;\frac{ q_i^+ - q_i^- } { q^+_{i,r} - q^-_{i,r} }.
\]

\end{definition}

Now we discuss the properties of the above two domain splitting schemes.

\begin{theorem} \label{t5}

Let $\{Q_r\}$ be a sequence of domains
generated as above by applying the maximal-cut scheme in each splitting.
Then,  we have that either $\{Q_r\}$ is a finite sequence, i.e.,
$\exists r_0 <
\infty$ such that the critical vertices of $Q_{r_0}$ differ
in no more than one coordinates, or
$\{Q_r\}$ is an infinite sequence such that
\[
\lim_{r \rightarrow \infty}||p(\Phi_{\cal D}(\delta),q^{\alpha_{r}})-
p(\Phi_{\cal D}(\delta),q^{\beta_{r}})||
 = 0
\]
and
\[
\lim_{r \rightarrow \infty} k_l(\delta,Q_r)=
\lim_{r \rightarrow \infty} k_u(\delta,Q_r)=
\lim_{r \rightarrow \infty} k_m(\delta,Q_r).
\]
Moreover, the same result follows if $\{Q_r\}$ is a sequence of domains
generated as above by applying the fair-cut scheme in each splitting.
\end{theorem}

\begin{pf} We only need to consider the case that $\{Q_r\}$
is an infinite sequence.
Decompose the coordinates index set
${\cal I}=\{1,\cdots,\ell\}$ as
${\cal I}={\cal I}_{f} \bigcup {\cal I}_{\infty}$ where
\[
{\cal I}_{f}=\{i \in {\cal I}:
\;[q_i^-,\;q_i^+]\;{\rm is\;divided\;finite\;many\;times} \}
\]
and
\[
{\cal I}_{\infty}=\{i \in {\cal I}:
\;[q_i^-,\;q_i^+]\;{\rm is\;divided\;infinite\;many\;times} \}.
\]
Obviously,
$\lim_{r \rightarrow \infty}||q^{\alpha_{r}}- q^{\beta_{r}}|| = 0$
for the case that ${\cal I}_{f}=\phi$.
We only need to consider the case that ${\cal I}_{f} \neq \phi,\;\;
{\cal I}_{\infty} \neq \phi$.
Note that $\exists {r}_{1} > 0$ such that $q_{i,r}^+=q_{i,{r}_{1}}^+,\;\;
q_{i,r}^-=q_{i,{r}_{1}}^-,\;\;\forall i \in {\cal I}_{f},\;\;
\forall r \geq {r}_{1}$.
Define
\[
\zeta=\min_{i \in {\cal I}_{f} } \; q_{i,r_1}^+ - q_{i,r_1}^-.
\]
Then
\[
\min_{i \in
{\cal I}_{f} } \; q_{i,r}^+ - q_{i,r}^- = \zeta > 0,\;\;\forall r \geq {r}_{1}.
\]
Note that $\exists r_2 > 0$ such that
\[
q_{i,r}^+ - q_{i,r}^- < \zeta,\;\;\forall i \in {\cal I}_{\infty},\;\;
\forall r \geq {r}_{2}.
\]
We claim that
${\cal U}_r \bigcap {\cal I}_{f}= \phi, \;\;\forall r > \max\{r_1,r_2\}$.
In fact, if this is not the case,
then $\exists i_\star \in {\cal U}_r \bigcap {\cal I}_{f}$ such that
\[
q_{i_\star,r}^+ - q_{i_\star,r}^- =
\max_{i \in {\cal U}_r} q_{i,r}^+ - q_{i,r}^- \geq \zeta
\]
because $q_{i,r}^+ - q_{i,r}^- < \zeta,\;\; \forall i
\in {\cal U}_r \bigcap {\cal I}_{\infty}$.  It follows that $Q_r$
is split as $Q_{r+1}$ and $Q_{r+1}^{'}$ by dividing interval
$[q_{i_\star,r}^-, \; q_{i_\star,r}^+]$, which contradicts to
$q_{i,r}^+=q_{i,{r}_{1}}^+,\;\;
q_{i,r}^-=q_{i,{r}_{1}}^-,\;\;\forall i \in {\cal I}_{f},\;\;
\forall r \geq {r}_{1}$.  Thus, the claim is true and it follows that
\[
||q^{\alpha_{r}}- q^{\beta_{r}}||^2=\sum_{i \in {\cal I}_{\infty}}
(q_i^{\alpha_{r}}- q_i^{\beta_{r}})^2 \leq \sum_{i \in {\cal I}_{\infty}}
(q_{i,r-1}^+ - q_{i,r-1}^-)^2,\;\;\forall r >\max\{r_1,r_2\}.
\]
Therefore, $\lim_{r \rightarrow \infty}||q^{\alpha_{r}}- q^{\beta_{r}}|| = 0$.
Since $p(z,q)$ is a continuous function of $q$,
it follows that
$\lim_{r \rightarrow \infty}||p(\Phi_{\cal D}(\delta),q^{\alpha_{r}})-
p(\Phi_{\cal D}(\delta),q^{\beta_{r}})||
 = 0$.  By the definition of $k_l$ and $k_u$, we have
 \[
\lim_{r \rightarrow \infty} k_l(\delta,Q_r)=
\lim_{r \rightarrow \infty} k_u(\delta,Q_r)=
\lim_{r \rightarrow \infty} k_m(\delta,Q_r).
\]

Similarly, to show that the same result follows if $\{Q_r\}$ is a
sequence of domains
generated as above by applying the fair-cut scheme in each splitting,
we only need to consider the case that ${\cal I}_{f} \neq \phi,\;\;
{\cal I}_{\infty} \neq \phi$.
Note that $\exists {r}_{3} > 0$ such that $q_{i,r}^+=q_{i,{r}_{3}}^+,\;\;
q_{i,r}^-=q_{i,{r}_{3}}^-,\;\;\forall i \in {\cal I}_{f},\;\;
\forall r \geq {r}_{3}$.  Define
\[
n_s=\max_{i \in {\cal I}_{f} } \;
\frac{q_i^+-q_i^-} {q_{i,r_3}^+ - q_{i,r_3}^-}.
\]
Then
\[
\max_{i \in {\cal I}_{f} } \;
\frac{q_i^+ - q_i^-} {q_{i,r}^+ - q_{i,r}^-} =n_s < \infty,
\;\;\forall r \geq {r}_{3}.
\]
Note that $\exists r_4 > 0$ such that
\[
\frac{q_i^+ - q_i^-} {q_{i,r}^+ - q_{i,r}^-}>
n_s,\;\;\forall i \in {\cal I}_{\infty},\;\;
\forall r \geq {r}_{4}.
\]
We claim that
${\cal U}_r \bigcap {\cal I}_{f}= \phi, \;\;\forall r > \max\{r_3,r_4\}$.
In fact, if this is not the case,
then $\exists i_\star \in {\cal U}_r \bigcap {\cal I}_{f}$ such that
\[
\frac{ q_{i_\star}^+ - q_{i_\star}^- }
{ q_{i_\star,r}^+ - q_{i_\star,r}^- } =
\min_{i \in {\cal U}_r}  \frac{q_i^+ - q_i^-}
{q_{i,r}^+ - q_{i,r}^-} \leq n_s
\]
because $\frac{ q_i^+ - q_i^- } { q_{i,r}^+ - q_{i,r}^- } > n_s,
\;\; \forall i \in {\cal U}_r \bigcap {\cal I}_{\infty}$.
It follows that $Q_r$
is split as $Q_{r+1}$ and $Q_{r+1}^{'}$ by dividing interval
$[q_{i_\star,r}^-, \; q_{i_\star,r}^+]$, which contradicts to
$q_{i,r}^+=q_{i,{r}_{3}}^+,\;\;
q_{i,r}^-=q_{i,{r}_{3}}^-,\;\;\forall i \in {\cal I}_{f},\;\;
\forall r \geq {r}_{3}$.  Thus, the claim is true and it follows that
\[
||q^{\alpha_{r}}- q^{\beta_{r}}||^2=\sum_{i \in {\cal I}_{\infty}}
(q_i^{\alpha_{r}}- q_i^{\beta_{r}})^2 \leq \sum_{i \in {\cal I}_{\infty}}
(q_{i,r-1}^+ - q_{i,r-1}^-)^2,\;\;\forall r >\max\{r_3,r_4\}.
\]
Therefore, by the same argument as in the maximal-cut schemes,
 the result follows.

\end{pf}

\begin{remark}
From the proof of the theorem we can see that both
domain splitting schemes guarantee
$||q^{\alpha_{r}}-q^{\beta_{r}}||
\rightarrow  0$ while allow $Q_r \rightarrow Q_\infty$
where $Q_\infty$ is not a single point in ${\bf R}^\ell$.  Clearly,
to make a subdomain converge to a single point requires much
more computational effort than to make
$||q^{\alpha_{r}}- q^{\beta_{r}}|| \rightarrow  0$.
As we can see later, $||q^{\alpha_{r}}- q^{\beta_{r}}||
\rightarrow  0$ leads to the existence of a
sequence of lower bounds (or upper bounds) converging to $k_m$.  Therefore,
an algorithm based on the maximal-cut (or fair-cut) splitting scheme will
reduce much computational effort in computing $k_m$ than other algorithms
based on making subdomains converge to a single point in ${\bf R}^\ell$.
From the proof, we can also see that the convergence will not follow
if the domain splitting is made along
the larger but not the largest edges of each subdomain.
It was remarked in \cite{G} that a fair-cut avoids the problem of getting into
very narrow and long subdomains which can decrease the convergence speed.
From the proof, we can see that it plays a role much more than affecting the
speed of convergence.  It is a sufficient
condition to the existence of a sequence of lower bounds
(or upper bounds) converging to $k_m$.
We would like to point out that the maximal-cut scheme has better worst case
convergence behavior than that of the fair-cut scheme.
\end{remark}

To see the efficiency of the maximal-cut (or the fair-cut) splitting scheme,
it is helpful to compare the image of the last subdomain resulted from the
the maximal-cut (or the fair-cut) splitting scheme and the image of the
last subdomain obtained by the finely subdivision.
The situation is shown in Figure~\ref{fig_333}.

\begin{figure}[htb]
\centerline{\psfig{figure=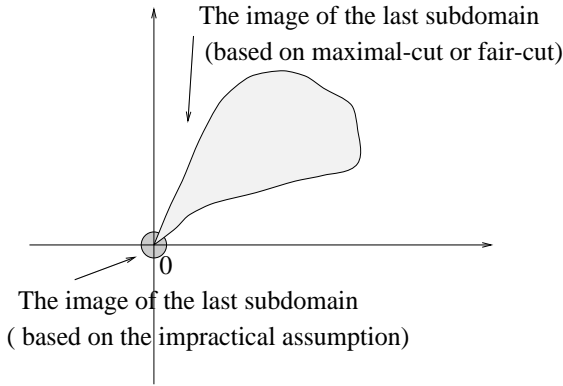,height=2in  ,width=3in  }}
\caption{The Image of the Last Subdomain}
 \label{fig_333}
\end{figure}

\section{Parallel Frequency Sweeping Algorithm}

In this section we shall investigate a new frequency sweeping structure.

\subsection{The Main Root of Inefficiency}

To the best of our knowledge, no effort in the existing literature has been
devoted to exploit a smart frequency sweeping strategy.
Existing techniques are basically as follows: Choose and grid a range of frequency.  Then calculate
 exactly $k_{m}$ for each gridded frequency.
 Finally, compare to find the minimum $k_{m}$ and
 return it as an estimate of $k_{max}$.

 For complicated root region ${\cal D}$,
  the number of frequencies to be evaluated for $k_m$ would be
  substantial in order to obtain a reasonably good estimate for $k_{max}$.
  Even the computation of $k_{m}$ for each frequency is very efficient,
 the overall complexity is still very high,
 because we need to evaluate $k_{m}$ for many frequencies.

 Thus for the sake of efficiency, it is natural to conceive
 a smart frequency sweeping strategy.  More specifically,
 we would raise the following question,

{\it Is it possible to obtain
the stability margin $k_{max}$
without tightly bounding $k_{m}(\delta_{j},{\cal Q})$
for each frequency $\delta_{j}$?}

The following section is devoted to answering this question.

\subsection{Parallel Frequency Sweeping Algorithm}

Consider the same
set of gridded frequencies $\delta_j,\;\;j=1,\cdots,n_rn_c$
defined by ~(\ref{grid}) and relabel them as
\[
\delta_{ij}:=
\delta_{l}+\frac{(\delta_{u}-\delta_{l})(i-1)}{n_r}+
 \frac{(\delta_{u}-\delta_{l})(j-1)}{n_rn_c}
,\;\;\;i=1,\cdots,n_r,\;\;j=1,\cdots,n_c.
\]

We are now in a position to present our Parallel Frequency Sweeping Algorithm
as follows.

{\bf Algorithm $2$ --- Parallel Frequency Sweeping Algorithm}

\begin{itemize}

\item Step $1$: Initialize. Set $j=1$. Set $\hat{k}=\infty$.
Set tolerance $\epsilon > 0$.  Set maximal iteration number $IT$.

\item Step $2$: Update $\hat{k}$ and record the number of
iterations $r(j)$ for frequency $\delta_{ij}$ by the following steps.

\begin{itemize}

\item Step $2$--$1$:
Let ${\cal U}_{ij}=\{Q_k\}={\cal Q},\;\; i=1,\cdots,n_r$. Set $r=1$.

\item Step $2$--$2$: If $r=IT+1$ or ${\cal U}_{ij}$ is empty for any
$i \in \{1,\cdots,n_r\}$ then record $r(j)=r$ and go to Step $3$,
else do the following for all $i$
such that ${\cal U}_{ij}$ is not empty.

\begin{itemize}

\item Choose $Q$ to be any element of ${\cal U}_{ij}$
with
\[
k_l(\delta_{ij},Q)=
\min_{Q_k \in {\cal U}_{ij}} \;k_l(\delta_{ij},Q_k).
\]

\item Partition $Q$ into $Q_a$ and $Q_b$ by applying a {\it maximal-cut}.

\item Remove $Q$ from ${\cal U}_{ij}$.

\item Update
\be
\hat{k}=\min \{\hat{k},\;
k_u(\delta_{ij},Q_a),\;k_u(\delta_{ij},Q_b)\}.
\label{update}
\ee

\item  Add any $Q \in \{Q_a, Q_b\}$ with two or more critical vertices
to ${\cal U}_{ij}$.

\item Remove from ${\cal U}_{ij}$ any $Q$ with
\be
0 \notin {\rm conv}
\left(p\left(\Phi_{\cal D}(\delta_{ij}),
\frac{\hat{k}}{1+\epsilon}Q\right) \right).
\label{con}
\ee

\end{itemize}

\item Step $2$--$3$: Set $r=r+1$ and go to Step $2$--$2$.

\end{itemize}

\item Step $3$: If $j=n_c$ then STOP, else set $j=j+1$ and go to Step $2$.

\end{itemize}

In Algorithm $2$, $n_r$ branches of frequency sweeping are performed
 in parallel with starting
frequencies $\delta_{i1},\;\;i=1,\cdots,n_r$ and step size
$\frac{\delta_{u}-\delta_{l}}{n_rn_c}$.
Each branch of frequency sweeping is not independent,
they exchange information.  The information is applied to determine
the subdomains to be eliminated from further consideration and to
update $\hat{k}$.
Finally, $\hat{k}$ is returned as the robust stability margin.
Algorithm $2$ is visualized
in the following Figure ~\ref{fig_18}.

\begin{figure}[htb]
\centerline{\psfig{figure=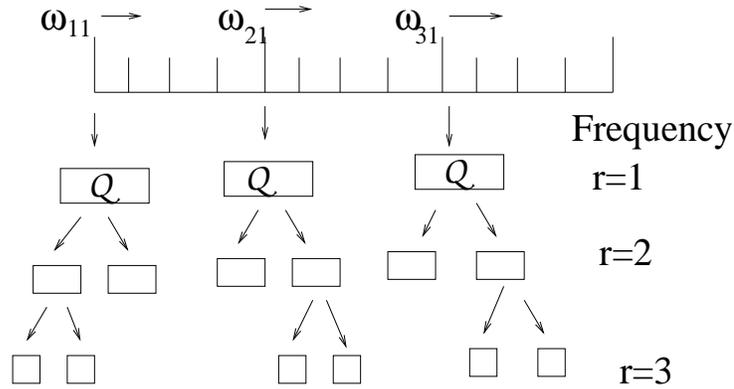,height=2.0in  ,width=3.8in  }}
\caption{A Picture of Parallel Frequency Sweeping Algorithm.  $N=3,\;\;K=4$.
 } \label{fig_18}
\end{figure}

\begin{remark}

As we can see from Step $2$--$2$,
there are two mechanisms which contribute to the efficiency of the Parallel
Frequency Sweeping Algorithm.  First, any subdomain $Q$ that satisfies condition (\ref{con}), which is equivalent to
\be
k_l(\delta_{ij}, Q) <
\frac{\hat{k}}{1+\epsilon},
\label{con8}
\ee
will never be partitioned again and thus can be
eliminated from further consideration.
Second, any subdomain $Q$ with critical vertices
differing in no more than one coordinates
will never be partitioned again and thus can be
eliminated from further consideration.

We would like to note that the proposed Parallel
Frequency Sweeping provides substantial improvement
on efficiency than the algorithms in \cite{GS}.
This can be explained by the significant relaxation of the
condition for eliminating a subdomain from consideration.
By ~(\ref{update}), ~(\ref{con}) and ~(\ref{con8}) we can see
that $\hat{k}$ is the minimum record of
the upper bounds among all subdomains evaluated (no matter it belongs to
the same frequency or not) and is contracted to
$\frac{\hat{k}}{1+\epsilon}$.  In contrast, in algorithms of \cite{GS}
the minimum record of upper
bounds is obtained for the frequency being considered only.  Therefore,
the condition for eliminating a subdomain from consideration
is much looser than its counterpart of algorithms in \cite{GS}.
Consequently,  such a significant relaxation results in
a substantial decrease
of the number of subdomains needed to be evaluated.

\end{remark}

\begin{remark}

In Algorithm $2$, at each iteration, only the domain with the
smallest lower bound is partitioned.  This mechanism differs from
that of Algorithm $1$ in which all domains are partitioned and thus the
number of domains increases rapidly.  We can see that
Algorithm $2$ effectively controls the growth of the number of subdomains and
thus is much efficient than the conventional algorithm.

\end{remark}

\begin{remark}
It is important to note that Algorithm $2$ involves only one CPU processor.
It is fundamentally different from the parallel algorithms
which involves more than one CPU processor.
\end{remark}

\begin{remark}
The speed of computation depends on the choice of integers $n_r$ and $n_c$. When the total number
of frequencies is fixed (i.e., $n_r n_c$ is constant),  the number of branches
of frequency sweeping $n_r$ should not be too small.
Small $n_r$ will hinder the improvement of efficiency benefited from the parallelism.
However, very large $n_r$ will not result in optimal performance either.
\end{remark}

In addition to the novel frequency sweeping strategy, another
character of Algorithm $2$ is that
there is no tolerance criteria directly forced on the final result,
however, the final result falls into tolerance automatically.

\begin{theorem} \label{t4}

Suppose that the maximum iteration number $IT=\infty$.
For arbitrary tolerance $\epsilon > 0$,
Algorithm $2$ stops with a finite number of
domain splittings for each $j$, i.e., $r(j) < \infty,\;\;\forall j$.
Moreover, the final $\hat{k}$ satisfies
\[
0 \leq \frac{\hat{k}- \tilde{k}_{max}}
{\tilde{k}_{max}}    < \epsilon.
\]
\end{theorem}

\begin{pf}
We first show the final $\hat{k} \geq \tilde{k}_{max}$.  Let $k_u$ be the
upper bound of domain $Q \subseteq {\cal Q}$ which ever
appeared during the execution of Algorithm $2$.
Let $\delta$ be the associated frequency of $Q$.
Note that
 $0 \in p(\Phi_{\cal D} (\delta), k_u Q) \subset
 p(\Phi_{\cal D} (\delta), k_u {\cal Q})$ and thus $k_u \geq \tilde{k}_{max}$.
 Note that the final $\hat{k}$ is the minimum record of all such $k_u$'s,
 thus $\hat{k} \geq \tilde{k}_{max}$.

 We next need to show that Algorithm $2$ stops with
a finite iteration number $r(j)$ for each $j$.
Suppose $\exists j$ such that $r$ goes to $\infty$.
Then $\exists \delta_{ij}$ such that $r$ goes to $\infty$.
Therefore, we can construct a sequence of nested domains $\{Q_{r}\}$ such that
$Q_1 \supset Q_2 \supset \cdots \supset Q_r \supset Q_{r+1} \supset \cdots $.  Thus by
Theorem ~\ref{t5} we have that $\forall \epsilon > 0,\;
\exists r_0 < \infty$ such that
 $k_{u}(\delta_{ij},Q_r)-k_{l}(\delta_{ij},Q_r) <
\frac{\epsilon}{1+\epsilon} \tilde{k}_{max},\;\;\forall r \geq r_0$.
Thus $\min\{\hat{k},k_{u}(\delta_{ij},Q_{r_0})\}-
k_{l}(\delta_{ij},Q_{r_0}) <
\frac{\epsilon}{1+\epsilon} \tilde{k}_{max}$.
Note that $k_{l}(\delta_{ij},Q_{r_0+1})
\geq k_{l}(\delta_{ij},Q_{r_0})$ because
 $Q_{r_0+1} \subset Q_{r_0}$.
 Also note that $\hat{k}$ never increases,
 thus we have $\hat{k} - k_{l}(\delta_{ij},Q_{r_0+1}) <
 \frac{\epsilon}{1+\epsilon} \tilde{k}_{max}$.
 Note that $\hat{k} \geq \tilde{k}_{max}$, we have
 \[
 \hat{k} - k_{l}(\delta_{ij},Q_{r_0+1}) <
 \frac{\epsilon}{1+\epsilon} \hat{k}\;\;\; \Longrightarrow \;\;\;
 \frac{\hat{k}}{1+\epsilon} < k_{l}(\delta_{ij},Q_{r_0+1}),
 \]
which implies that $0 \notin {\rm conv} \left(p \left(\Phi_{\cal D}(\delta_{ij}),
\frac{\hat{k}}{1+\epsilon}Q_{r_0+1} \right) \right)$.
Therefore, by Algorithm $2$
$Q_{r_0+1}$ will not be splitted.  This is a contradiction.
 So, we have shown that Algorithm $2$ stops with a finite number of
domain splittings for each $j$.

Note that $\exists \delta_{ij}$ such that
$k_{m}(\delta_{ij},{\cal Q})=\tilde{k}_{max}$.
Moreover, $\exists q^\star \in {\cal Q}$ such that
$p(\Phi_{\cal D}(\delta_{ij}), \tilde{k}_{max} q^\star)=0$.
Since Algorithm $2$ stops with a finite number of
domain splittings for each $j$, we have that all subdomains
ever generated are finally eliminated from consideration.  Thus,
there must exists $Q^{\star}$ which
contains $q^\star$ be eliminated from consideration at
a certain level of splitting.

Assume that $\hat{k}= \bar{k}$ when $Q^{\star}$ is eliminated
from consideration.
Then either
$0 \notin p(\Phi_{\cal D}(\delta_{ij}), \frac{\bar{k}}{1+\epsilon}Q^{\star})$
or the critical vertices of $Q^{\star}$
differ in no more than one coordinates.
If the first case is true, then by
$p(\Phi_{\cal D}(\delta_{ij}), \tilde{k}_{max} q^\star)=0$ and $q^\star \in Q^{\star}$,
 we have that
$0 \in p(\Phi_{\cal D}(\delta_{ij}), \tilde{k}_{max} Q^{\star})$.
Thus by $0 \notin p(\Phi_{\cal D}(\delta_{ij}),\frac{\bar{k}}{1+\epsilon} Q^{\star})$,
  we have $\frac{\bar{k}}{1+\epsilon}< \tilde{k}_{max}$.  Obviously,
the final $\hat{k} \leq \bar{k}$ and thus
\[
\frac{\hat{k}}{1+\epsilon}< \tilde{k}_{max} \;\;\;
\Longrightarrow \;\;\;\frac{\hat{k}- \tilde{k}_{max}}
{\tilde{k}_{max}}    < \epsilon.
\]
If the latter case is true, then we have
\[
k_l(\delta_{ij}, Q^{\star})=
k_u(\delta_{ij}, Q^{\star})=\tilde{k}_{max}=\hat{k}.
\]

The proof is thus completed.

\end{pf}

\begin{remark}
From the proof, we can see that the existence of a sequence of
upper bounds converging to $\tilde{k}_{max}$
is due to the convergence of the distance of critical vertices.
\end{remark}

\begin{remark}
By specifying $\epsilon$ in Algorithm $2$, we can obtain an estimate $\widehat{k}$ such that $
0 \leq \frac{\hat{k}- \tilde{k}_{max}}
{\tilde{k}_{max}}    < \epsilon$.
In this sense, we can say that Algorithm 2 provides a global solution for
searching $\tilde{k}_{max}$.  However, it should be noted that it is not
necessary a global solution for the exact robust ${\cal D}$-stability margin $k_{max}$.
This is because it is impossible to search
the whole range of the generalized frequency $\delta$.
It is only feasible to perform the search over a set of discrete values of $\delta$.
\end{remark}

\begin{theorem} \label{t6}

Suppose that the maximum iteration number $IT<\infty$ and that
Algorithm $2$ stops with $r(j) < IT,\;\;\forall j$.  Then
the final $\hat{k}$ satisfies
\[
0 \leq \frac{\hat{k}- \tilde{k}_{max}}
{\tilde{k}_{max}}    < \epsilon.
\]
\end{theorem}
\begin{pf} Since Algorithm $2$ stops with $r(j) < IT,\;\;\forall j$,
we can conclude that all subdomains ever generated are finally
eliminated from consideration.
Thus the result follows from similar argument as for Theorem ~\ref{t4}.
\end{pf}

\section{An Illustrative Example}

Our computational experience shows that Algorithm $2$
provides a significant improvement upon conventional algorithms for
 most control problems.
Moreover, the improvement depends on the problems and
can be arbitrarily good.  To illustrate, we consider
an example with ${\cal D}$ chosen
to be the open left half plane.
The applications to the problems with
complicated ${\cal D}$ are in exactly the same spirit.

The state space equation of the linear system is as follows:
\[
\begin{array}{l}
\dot{x} = A(q) x + B u\\
y= C x
\end{array}
\]
where
\[
A(q)=\left[\begin{array}{cccc}
-1-0.5q_1 & -10 & -1 & 10\\
-0.5 & -1 & 1 & 0.5\\
0.5 & -4 & -1 & -10\\
-10 & 0.5 & 0 & -2.5-1.5q_2\end{array}\right],\;\;\;B=\left[\begin{array}{c}
1 \\
2 \\
1 \\
1 \end{array}\right],\;\;\;C=\left[\begin{array}{cccc}
1 & 0.5 & 1 & 0.5
\end{array}\right]
\]
with uncertain parameter $q \in
{\cal Q}=[0,1] \times [0,3] \subset {\bf R}^2$.
We obtained a polynomial for this system as
$p(s,q)=det(sI-A(q))$ which has a multilinear structure.

The upper bound and lower bound of $k_m$ on ${\cal Q}$
is shown in
Figures $4-7$.  We can see that for most of the frequencies the
 upper bounds and
 lower bounds are far apart and thus the importance of domain splitting
  is obvious.

\begin{figure}[htb]
\centerline{\psfig{figure=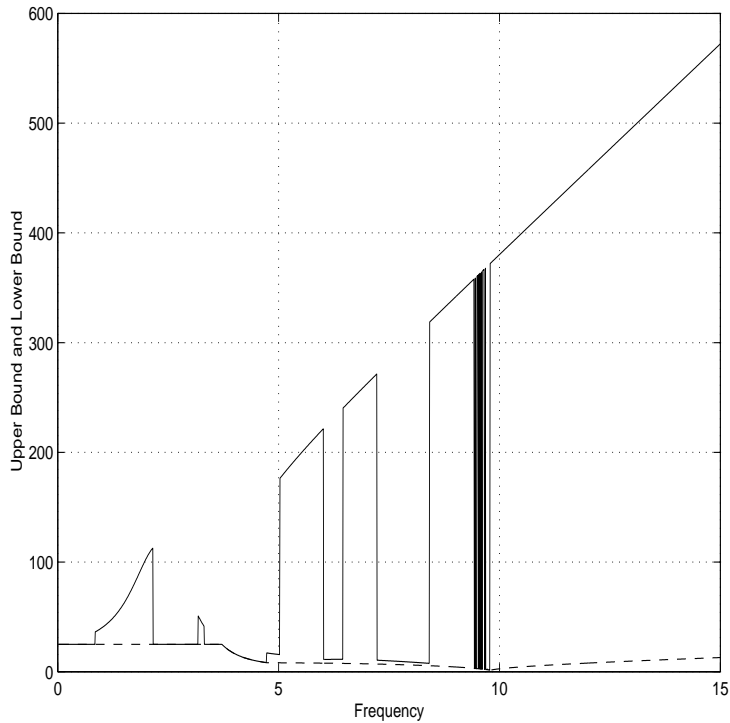,height=3.8in  ,width=3.8in  }}
\caption{
$k_m$ upper bound and lower bound on ${\cal Q}$.
The upper bound is plotted in dashed
line and the lower bound is plotted in solid line. } \label{fig_5}
\end{figure}

\begin{figure}[htb]
\centerline{\psfig{figure=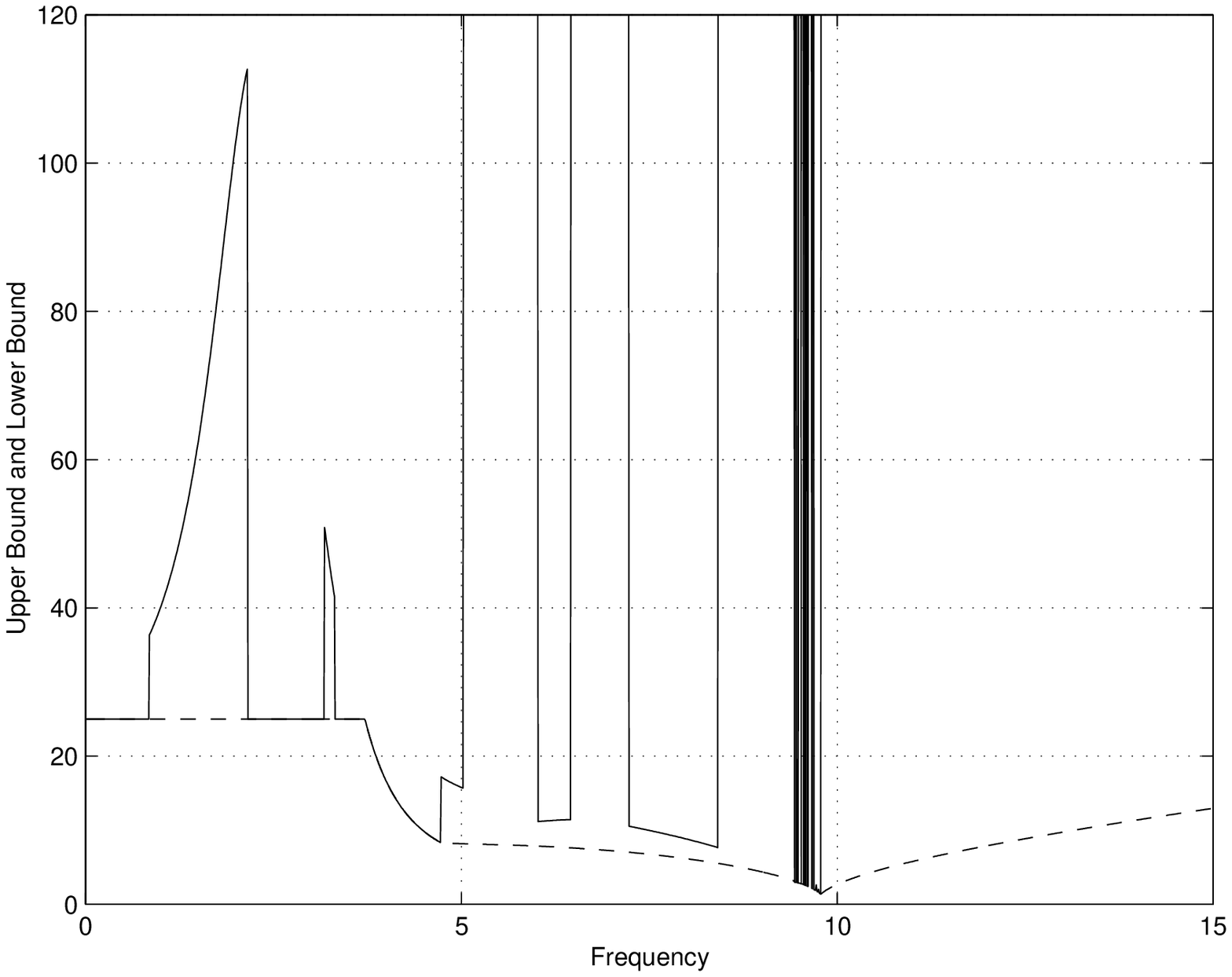,height=3.8in  ,width=3.8in  }}
\caption{
$k_m$ upper bound and lower bound on ${\cal Q}$.
The upper bound is plotted in dashed
line and the lower bound is plotted in solid line. } \label{fig_1}
\end{figure}

\begin{figure}[htb]
\centerline{\psfig{figure=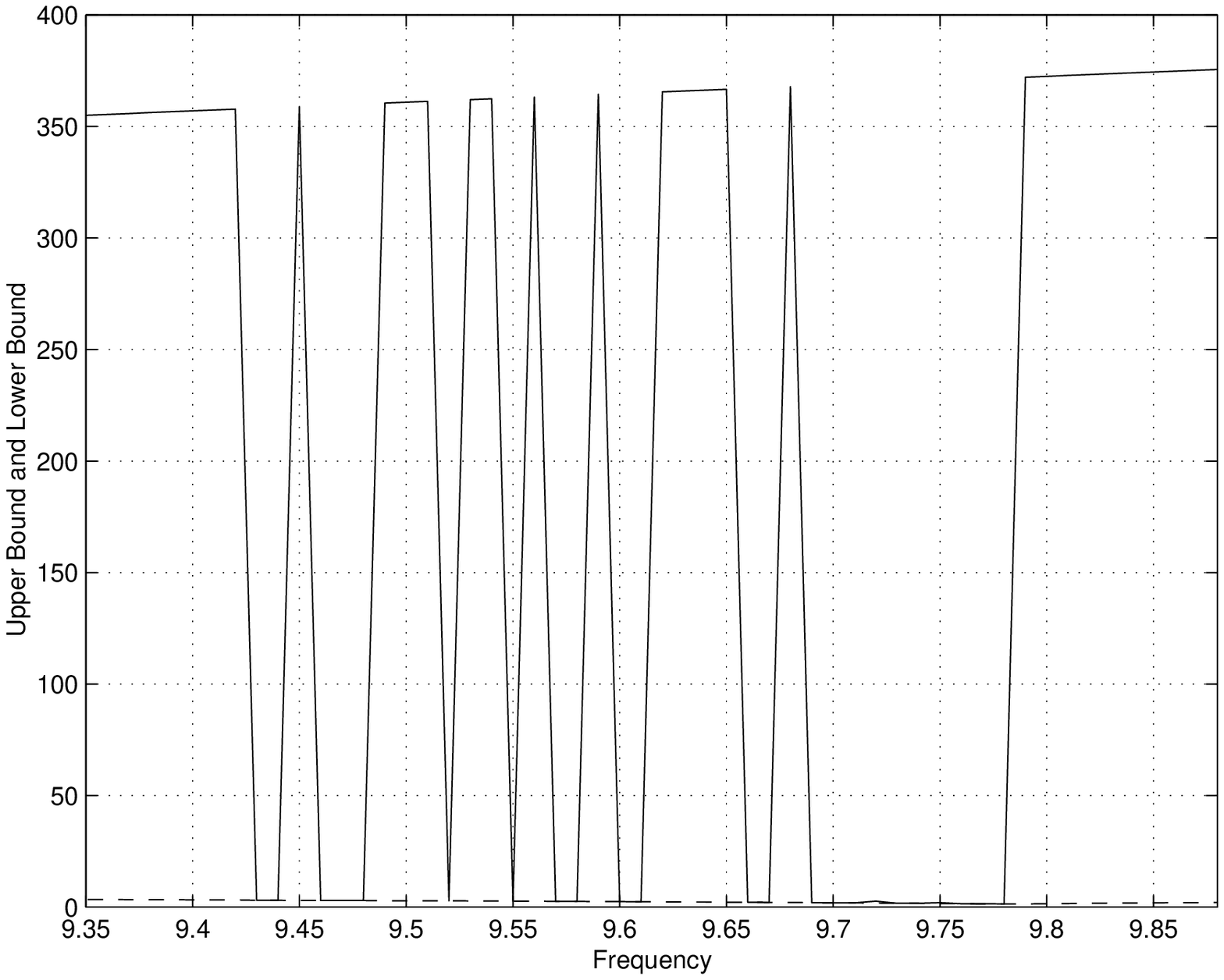,height=3.8in  ,width=3.8in  }}
\caption{
$k_m$ upper bound and lower bound on ${\cal Q}$.
The upper bound is plotted in dashed
line and the lower bound is plotted in solid line. }
\label{fig_2}
\end{figure}

\begin{figure}[htb]
\centerline{\psfig{figure=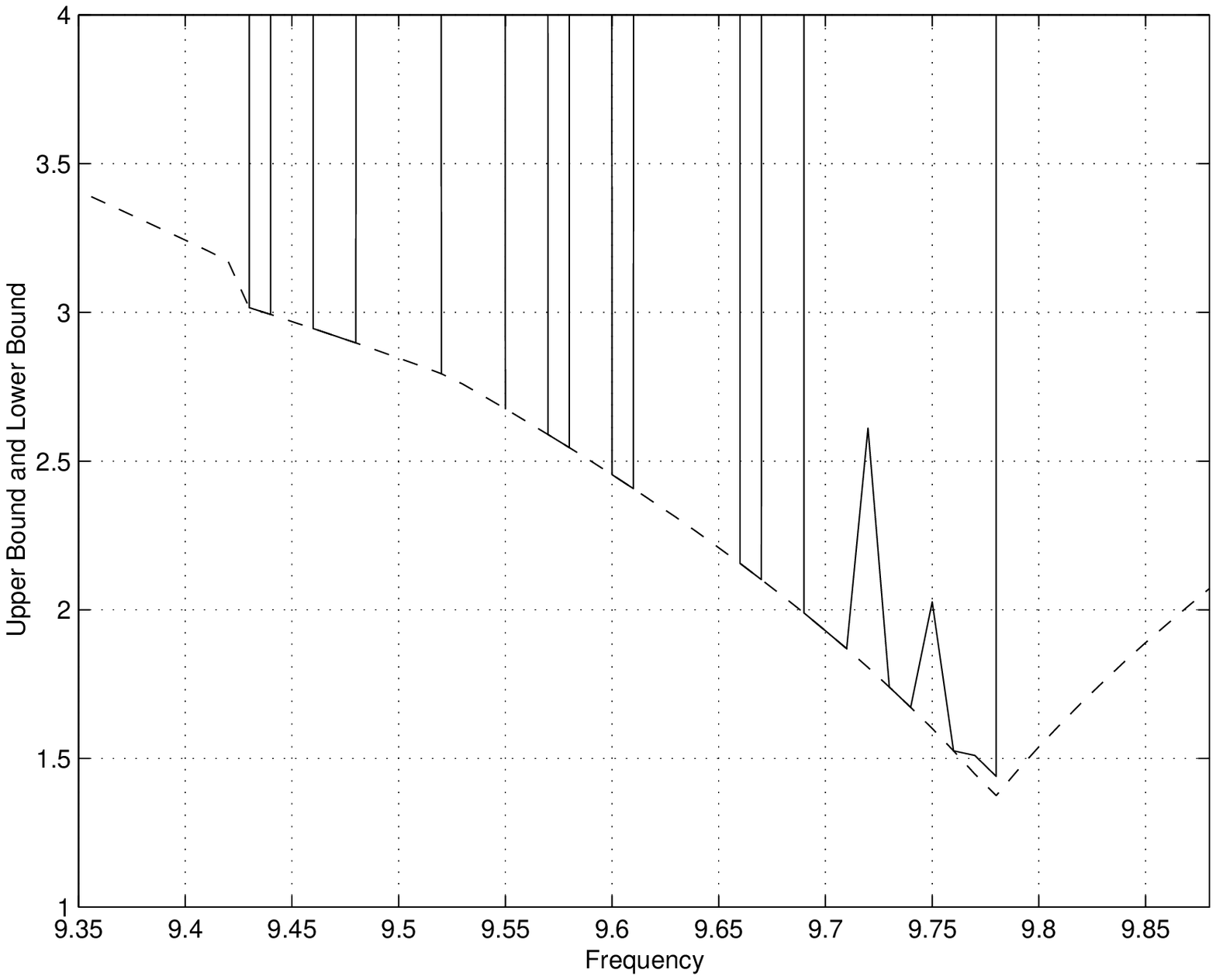,height=3.8in  ,width=3.8in  }}
\caption{
$k_m$ upper bound and lower bound on ${\cal Q}$.
The upper bound is plotted in dashed
line and the lower bound is plotted in solid line. }
\label{fig_4}
\end{figure}

To compute $k_{max}$, we uniformly grid frequency
band $[0.01 ,\;15.01]$ and obtain $1,500$ grid frequencies as
\[
\omega_j=0.01j, \;\;j=1,\cdots,1500.
\]
In Algorithm $2$, we choose the relative error $\epsilon=0.01$
and $n_r=30,\;\;n_c=50$.
The $1,500$ frequencies are regrouped as
\[
\omega_{ij}=0.01+0.2(i-1)+0.01(j-1),\;\;i=1,\cdots,30;\;\;j=1,\cdots,50.
\]
We ran the program in a Sun Spark work-station.  The running time
is about $80$ seconds.
The total number of domains evaluated is $1,570$.
We obtained $\hat{k}=1.4384$ which is achieved
at frequency $\omega_{20,18}=9.78$.  By Theorem ~\ref{t4},
we can concluded that
\[
0 \leq \frac{\hat{k}-\tilde{k}_{max}}
{\tilde{k}_{max}}    < \epsilon =0.01 .
\]

To compare the performance of the conventional algorithm with
that of Algorithm $2$,
it is fair to choose the tolerance $\varepsilon_r=0.01$ in Algorithm $1$.
We also
ran the program in the same work-station.  The running time
is about $9$ hours.
The total number of domains evaluated is $64,813$.
We obtained $\hat{k}=1.4380$ which is achieved
at frequency $\omega_{9783}=9.783$.  Therefore,
Algorithm $2$ has a speed-up of $400$ over the conventional algorithm.
Moreover, the number of
domains evaluated in Algorithm $2$ is only a small fraction
(which is $\frac{1570}{64813} \approx 0.0242$) of that of the
 conventional algorithm.  The number of domains evaluated in
 Algorithm $2$ and the conventional
 one for each frequency is shown respectively
 in Figure~\ref{fig_11} and Figure~\ref{fig_6}.

\begin{figure}[htb]
\centerline{\psfig{figure=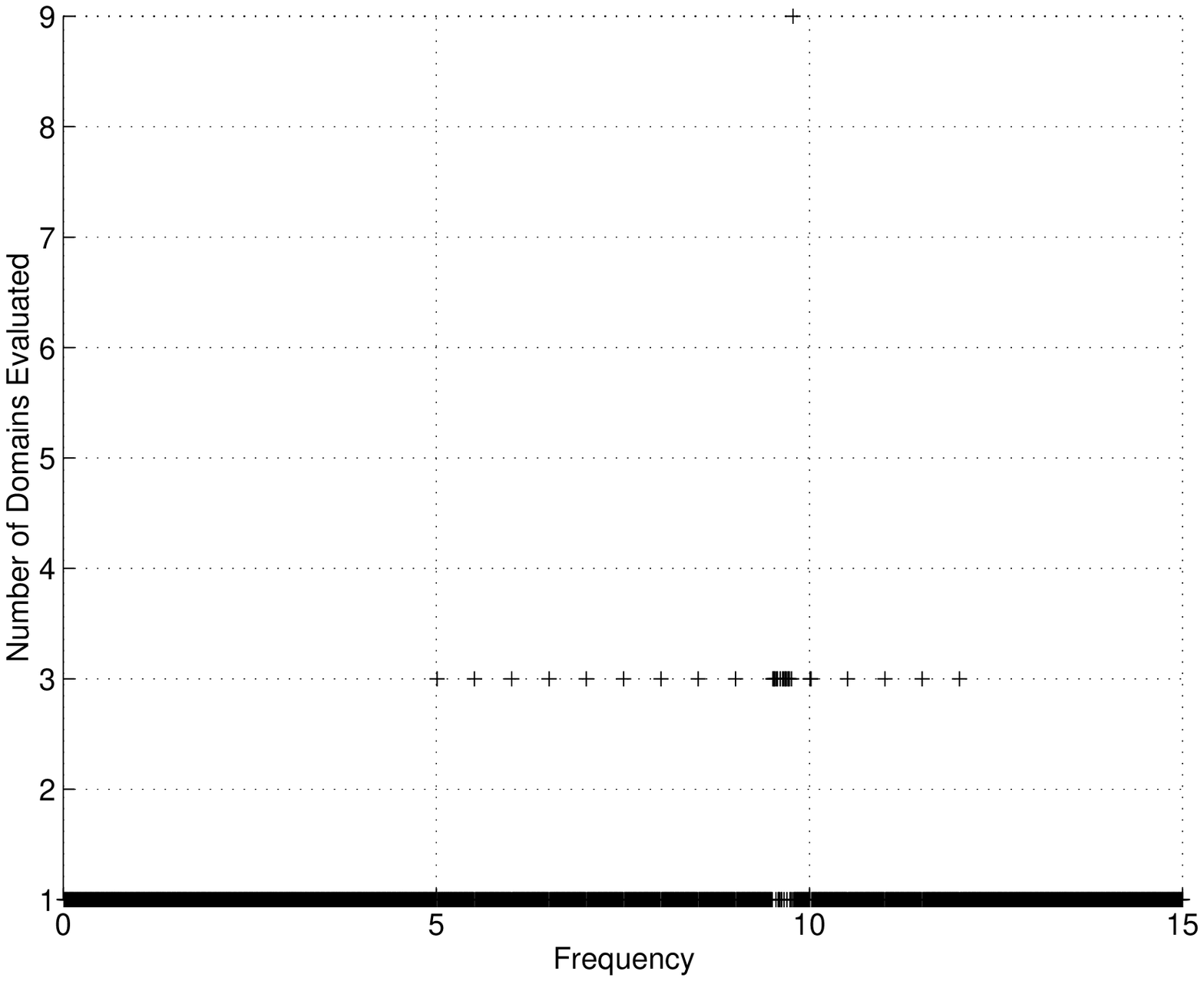,height=3.8in  ,width=3.8in  }}
\caption{Domains evaluated in Algorithm $2$.}
 \label{fig_11}
\end{figure}

\begin{figure}[htb]
\centerline{\psfig{figure=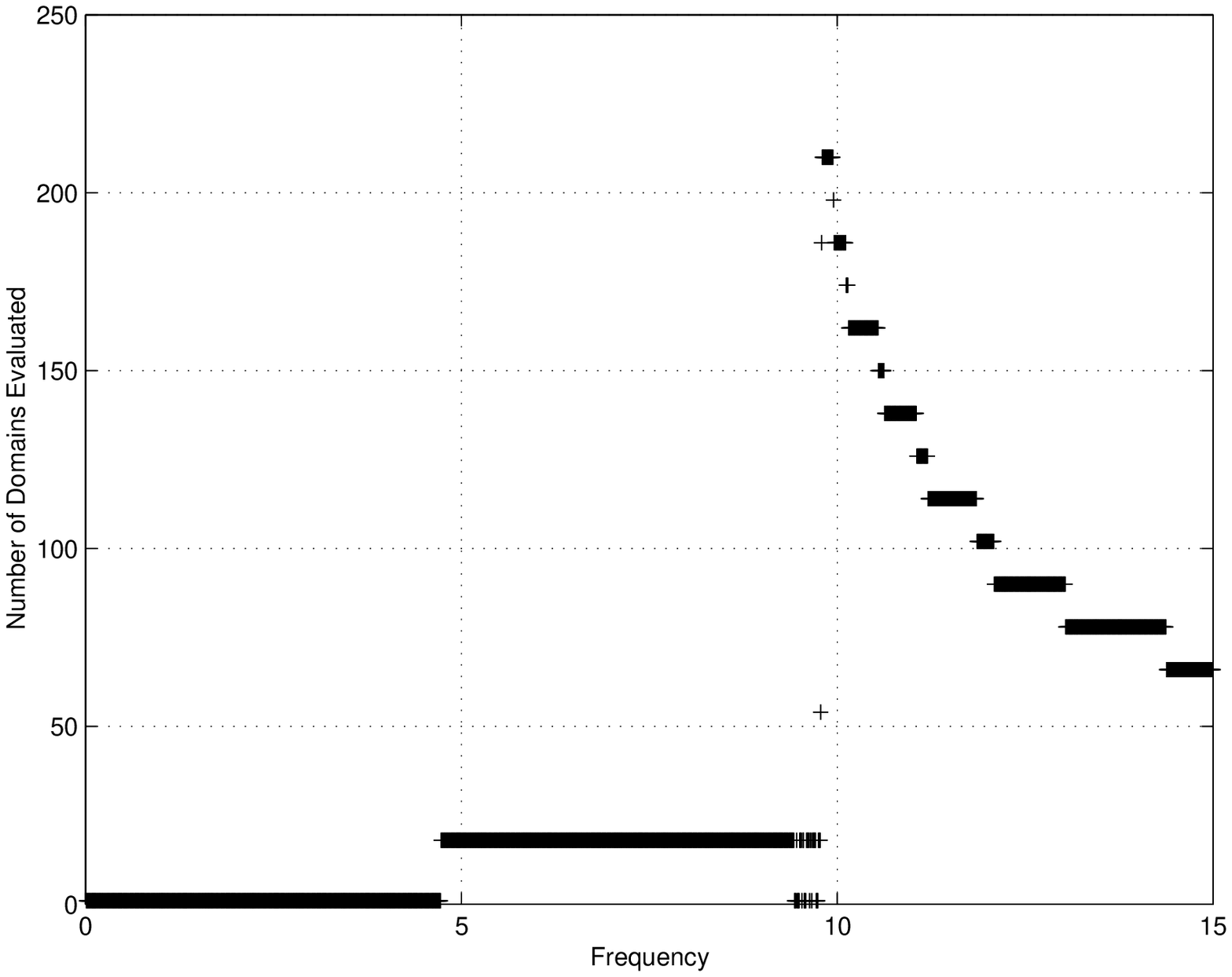,height=3.8in  ,width=3.8in  }}
\caption{Domains evaluated in the conventional algorithm.} \label{fig_6}
\end{figure}

 We can see that Algorithm $2$ provides much superior
  performance than the conventional algorithms.  The improvement comes
  from the characteristic eliminating mechanisms in Algorithm $2$.
  More formally, we describe the eliminating process
  in Algorithm $2$ as follows.

Let $U$ be a record of the global
upper bound achieved by frequency $\omega$.  Let
$Q_{ij} \subseteq {\cal Q}$ be a domain associated with
frequency $\omega_{ij}$.
When $Q_{ij}$ is eliminated, i.e.,
$k_l( \omega_{ij}, Q_{ij}) > \frac{U}{1+\epsilon}$ is satisfied,
there are only three cases as follows.

\begin{itemize}

\item Case (i): $\omega_{ij} < \omega$.  We call the elimination as {\it Backward Pruning}.

\item Case (ii): $\omega_{ij} > \omega$.  We call the elimination as {\it Forward Pruning}.

\item Case (iii): $\omega_{ij} = \omega$.  We call the elimination as {\it Present Pruning}.

\end{itemize}

All the above three types of pruning processes play important roles in
Algorithm $2$.  However, there is only {\it Present Pruning} in
the conventional algorithm.  Therefore, Algorithm $2$ has
a much powerful pruning mechanism and is much more efficient.

In this example, we have $24$ records which are shown in
 Figure ~\ref{fig_10}.  The effectiveness of the
 three types of pruning processes are shown
 respectively in Figures~\ref{fig_9},~\ref{fig_7} and ~\ref{fig_8}.

\begin{figure}[htb]
\centerline{\psfig{figure=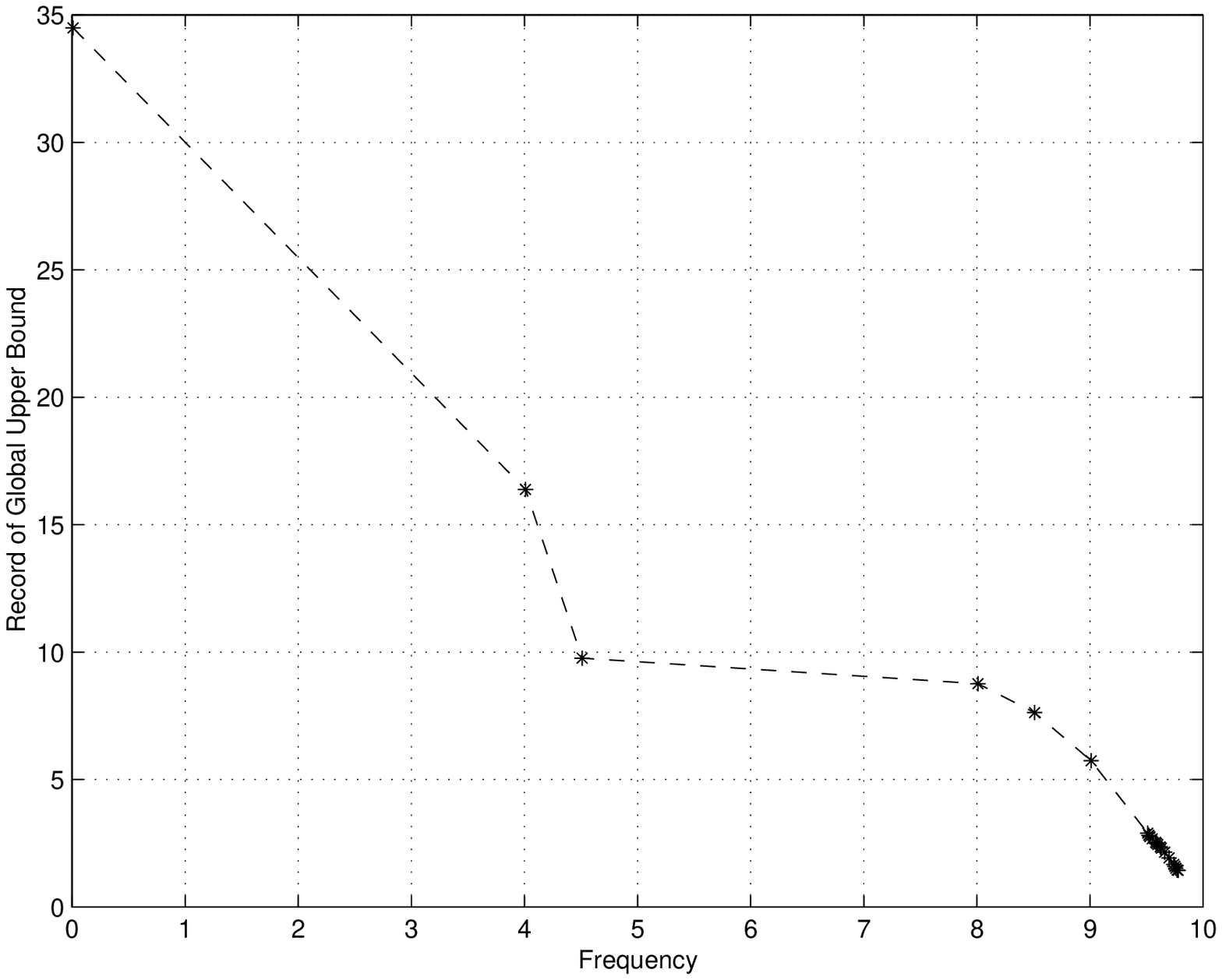,height=3.8in  ,width=3.8in  }}
\caption{Evolution of the global upper bound.
The $y$-coordinate represents the record value and the
 $x$-coordinate represents the frequency
achieving it.  Two consecutive records are connected by dashed line.
} \label{fig_10}
\end{figure}

\begin{figure}[htb]
\centerline{\psfig{figure=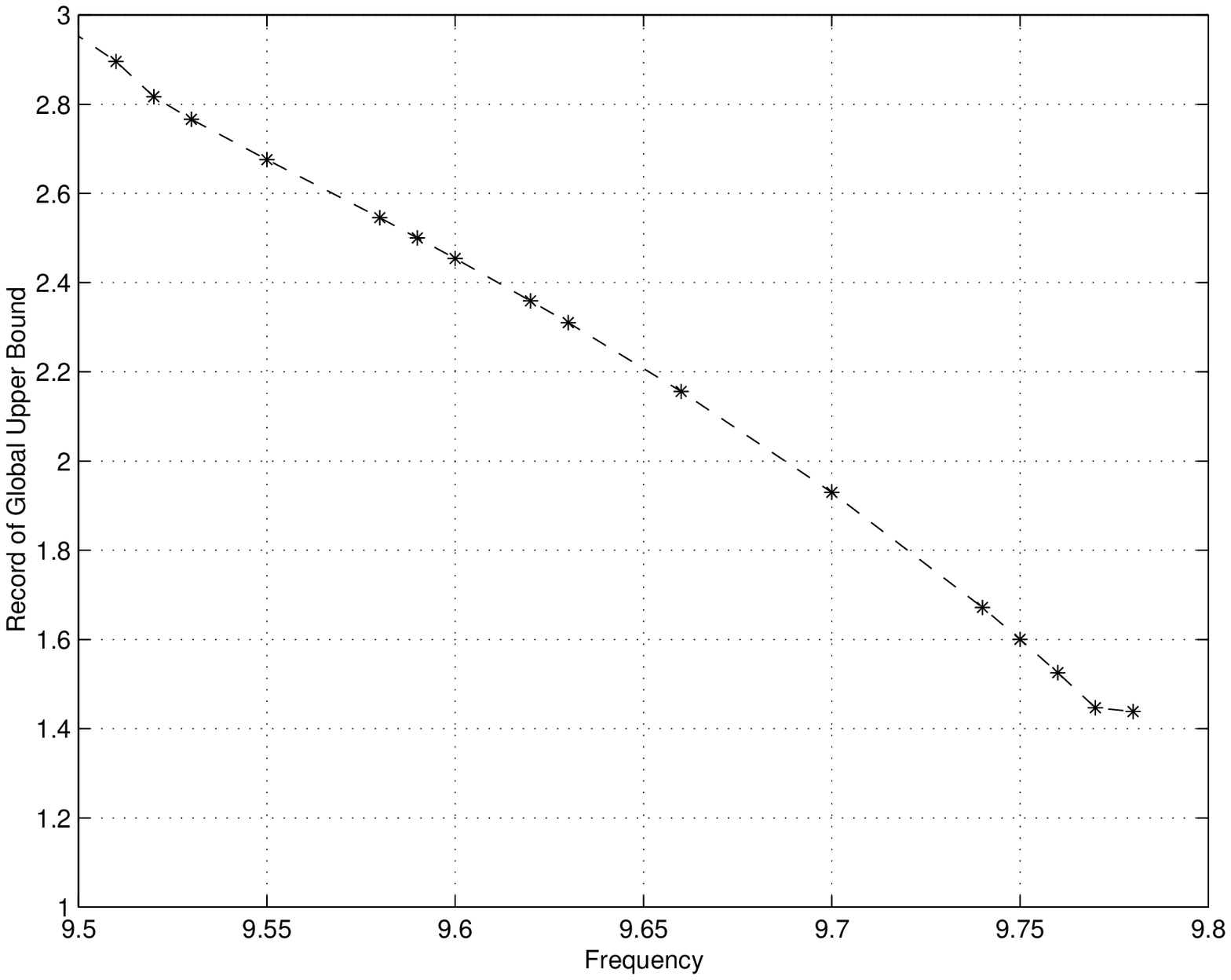,height=3.8in  ,width=3.8in  }}
\caption{Evolution of the global upper bound.
The $y$-coordinate represents the record value and the
 $x$-coordinate represents the frequency
achieving it.  Two consequent records are connected by dashed line.
} \label{fig_100}
\end{figure}

\begin{figure}[htb]
\centerline{\psfig{figure=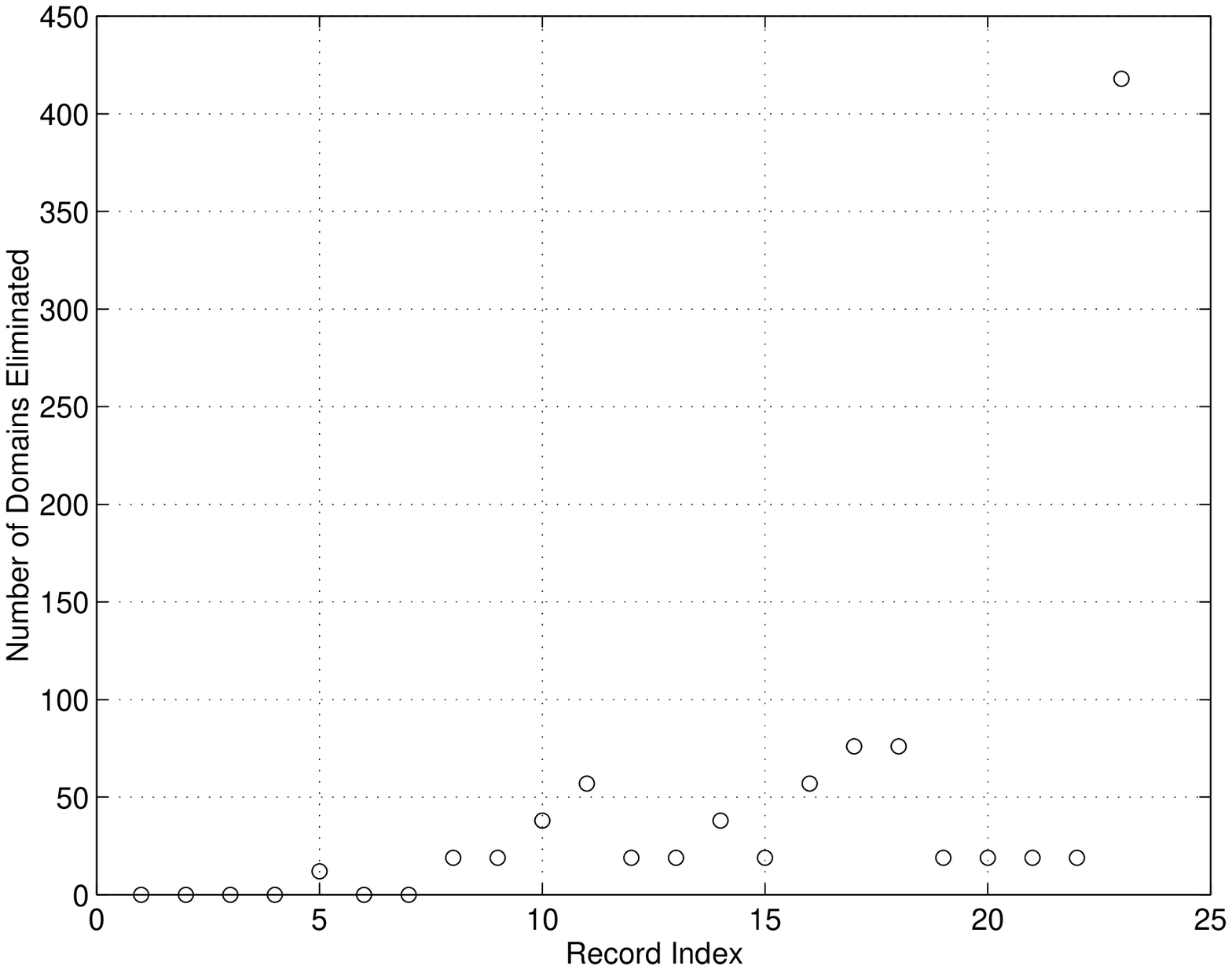,height=3.8in  ,width=3.8in  }}
\caption{Backward Pruning.  The $y$-coordinates represents
the number of domains eliminated by the record as Case (i).
The $x$-coordinate represents the record index.
} \label{fig_9}
\end{figure}

\begin{figure}[htb]
\centerline{\psfig{figure=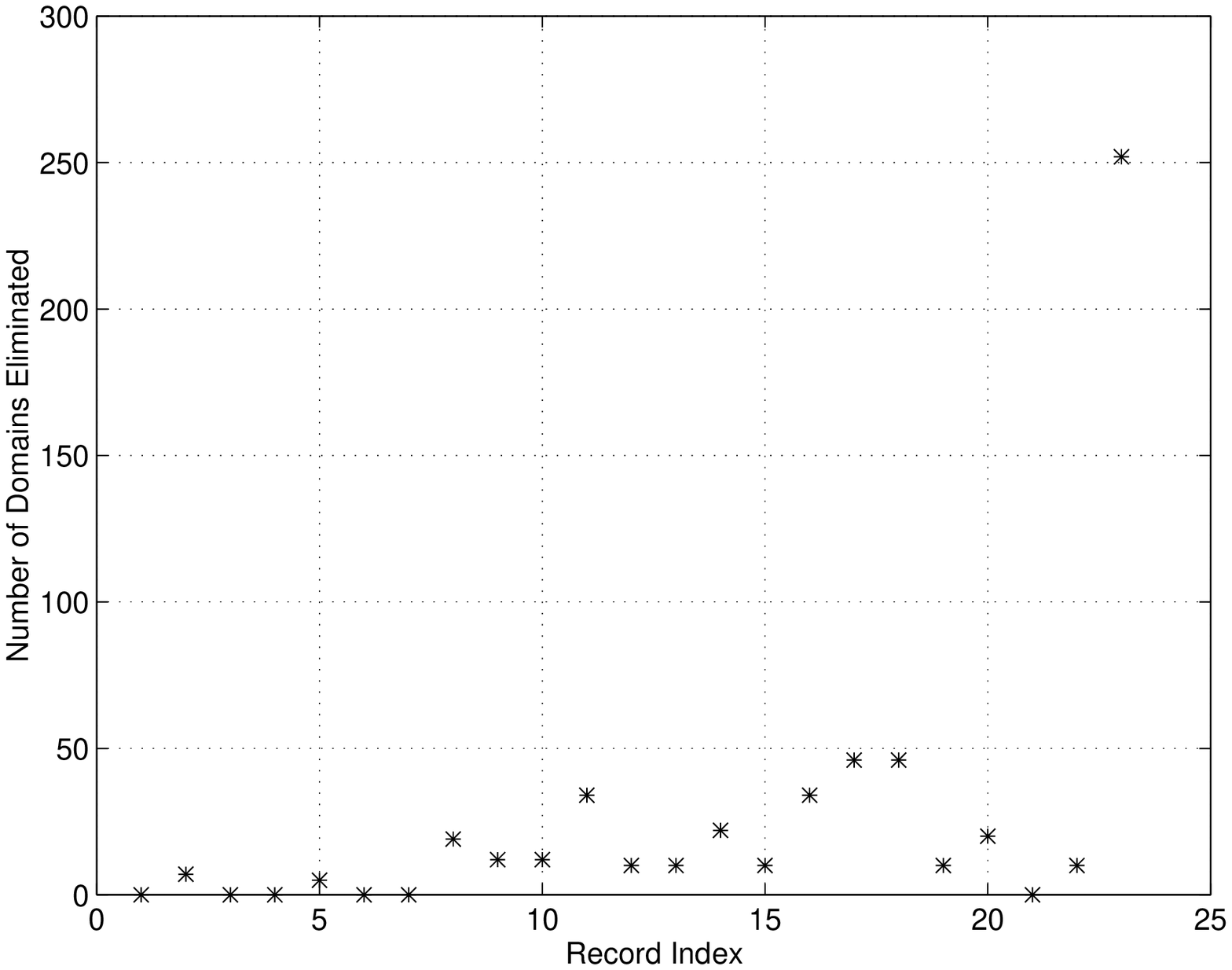,height=3.8in  ,width=3.8in  }}
\caption{Forward Pruning. The $y$-coordinates represents
the number of domains eliminated by the record as Case (ii).
The $x$-coordinate represents the record index.} \label{fig_7}
\end{figure}

\begin{figure}[htb]
\centerline{\psfig{figure=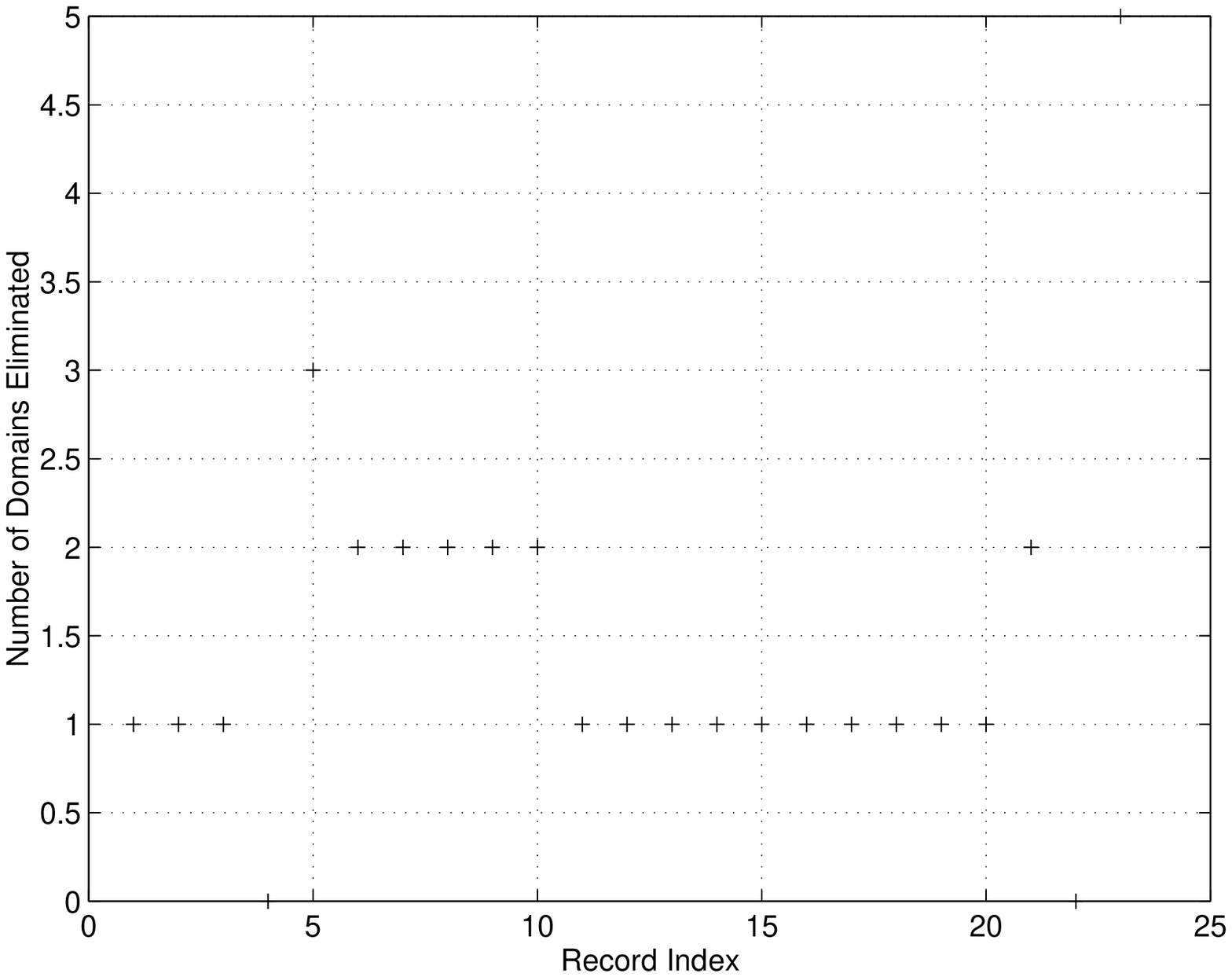,height=3.8in  ,width=3.8in  }}
\caption{Present Pruning.  The $y$-coordinates represents
the number of domains eliminated by the record as Case (iii).
The $x$-coordinate represents the record index.} \label{fig_8}
\end{figure}

\section{Conclusion}

We have developed techniques such as,
a parallel frequency sweeping strategy,
different domain splitting schemes,
which significantly reduce the computational complexity
and guarantee the convergence.  Our computational experience shows that
Algorithm $2$ provides a substantial improvement on efficiency in
comparison with the conventional algorithms.


\begin{thebibliography}{10}

\bibitem{Bal} Balas, G., Doyle, J. C., Glover, K., Packard, A.
 and Smith, R., {\it $\mu$-Analysis and Synthesis Toolbox},
 MUSYN Inc. and The MathWorks, Inc. 1995.

\bibitem{Ba} Balas, G. and Packard A.,
``The structured singular value $(\mu)$ framework,''
{\it The Control Hand Book}, pp. 671-687, CRC Press,
Inc, 1996.

\bibitem{B} Barmish, B. R.,  {\it New Tools for Robustness of Linear Systems},
Macmillan, New York, 1994.

\bibitem{BKT} Barmish, B. R., Khargonekar P. P., Shi Z. C. and Tempo R.,
``Robustness margin needs not be a continuous function of the problem data'',
{\it System and Control Letters}, Vol. 15, 1990.

\bibitem{BYDM} Braatz R. D., Young P. M., Doyle J. C. and Morari M.,
``Computational complexity of $\mu$ calculation'',
{\it IEEE Transactions on Automatic Control}, Vol. 39, No. 5, 1994.


\bibitem{C} Chen, C. T., {\it Linear System Theory and Design},
Holt, Rinehart and Winston, New York, 1984.

\bibitem{D} Doyle, J. C., ``Analysis of feedback system with structured
uncertainty,''
{\em IEE Proc.}, pt. D, vol. 129, no. 6, pp. 242-250, 1982.

\bibitem{FD} Frazer, R. A. and Duncan, W. J., ``On the criteria
for the stability of small motion,''
{\em Proceedings of the Royal Society},
A, vol. 124, pp. 642-654, London, England, 1929.


\bibitem{GS} de Gaston, R. R and M. G. Safonov, M. G.,
``Exact calculation of the multiloop stability margin,''
{\em IEEE Trans. Autom. Control}, vol. 33, pp. 156-171, 1988.

\bibitem{G} de Gaston, R. R.,
{\it Nonconservative Calculation of the Multiloop Stability Margin},
Ph.D. Dissertation,
University of Southern California, 1985.

\bibitem{K} Kharitonov, V. L., ``Asymptotic stability of an
equilibrium position of a family of systems of
linear differential equations,''
{\it Differentsial'nye Uravneniya}, vol. 14, pp. 2068-2088, 1978a.

\bibitem{NY} Newlin M. P. and Young P. M.,
``Mixed $\mu$ problems and branch and bound techniques,''
{\it Proc. IEEE Conference on Decision and Control
,} pp. 3175-3180, Tucson, Arizona, 1992.

\bibitem{Rohn} Rohn J. and Poljak R., ``Checking robust nonsigularity is NP hard'',
{\it Mathematics of Control, Signals and Systems}, 1992.


\bibitem{Pena} Sanchez Pena, R. S., {\it Robust Analysis of Feedback Systems with Parametric
and Dynamic Structured Uncertainty}, PhD Thesis, Caltech, 1988.

\bibitem{PS} Sanchez Pena, R. S. and Sideris, A., ``A general program to compute the
multivariable stability margin for systems with parametric uncertainty,''
{\it Proc. American Control Conference}, pp. 317-322, Atlanta, GA, 1988.


\bibitem{SG} Sideris, A. and de Gaston, R. R., ``Multivariable stability
margin calculation with uncertain correlated parameters,''
{\it Proc. IEEE Conference on Decision and Control
}, pp. 766-771, Athens, Greece, 1986.


\bibitem{S} Sondergeld, K. P., `` A generalization of the Routh-Hurwitz
criteria and an application to a problem in robust controller design,''
{\em IEEE Trans. Autom. Control}, vol. 28, pp. 965-970, 1983.

\bibitem{SP} Sideris, A. and Sanchez Pena, R. S., ``Fast calculation of the
multivariable stability
margin for real interrelated uncertain parameters,''
{\em IEEE Trans. Autom. Control}, vol, 34, pp. 1272-1276, 1989.


\bibitem{TB} Tempo, R. and Blanchini, F., ``Robustness analysis with real
parametric uncertainty,'' {\em The Control Hand Book}, pp. 495-505, CRC Press,
Inc, 1996.

\bibitem{VTM} Vicino, A., Tesi, A. and Milanese, M., ``Computation of
nonconservative stability perturbation bounds for systems
with nonlinearly correlated uncertainties,''
{\em IEEE Trans. Autom. Control}, 35, 835-841, 1990.

\bibitem{ZD} Zadeh, L. A. and Desoer, C. A.,
{\it Linear System Theory}, New York: McGraw-Hill, 1963.

\bibitem{ZDG} Zhou, K., Doyle J. C., and Glover K.,
{\it Robust and Optimal Control}, Prentice Hall,
Upper Saddle River, NJ, 1996.


\end{thebibliography}
\end{document}